\documentclass[final]{siamltex}

\usepackage{moreverb,rotating,graphics}
\usepackage{amsmath,amssymb,amsfonts}

\newtheorem{remark}{Remark}

\def\RR{{\mathbb R}}

\title{Stable parareal in time method for first and second order hyperbolic system\thanks{This
        work was supported by  ANR-06-CIS6-007.}}

\author{Xiaoying Dai \thanks{LSEC, Institute of Computational Mathematics and
  Scientific/Engineering Computing,
 Academy of Mathematics and Systems
  Science, Chinese Academy of Sciences, Beijing 100190, China ({\tt daixy@lsec.cc.ac.cn}).}
        \and Yvon Maday\footnotemark[3]\ \footnotemark[4]\ \footnotemark[5]}

\begin{document}

\maketitle

\renewcommand{\thefootnote}{\fnsymbol{footnote}}

\footnotetext[3]{UPMC Univ Paris 06, UMR 7598, Laboratoire Jacques-Louis Lions, F-75005, Paris, France}
\footnotetext[4]{CNRS, UMR 7598, Laboratoire Jacques-Louis Lions, F-75005, Paris, France}
\footnotetext[5]{Division of Applied Mathematics, Brown University,  Providence, RI, USA}

\begin{abstract}
The parareal in time algorithm allows to perform parallel simulations of time dependent problems. This algorithm has been implemented on many types of time dependent problems with some success. Recent contributions have allowed to extend the domain of application of the parareal in time algorithm so as to handle long time simulations of Hamiltonian systems. This improvement has managed to avoid the fatal large lack of accuracy of the plain parareal in time algorithm consequence of the fact that the plain parareal in time algorithm does not conserve invariants. A somehow similar difficulty occurs for problems where the solution lacks regularity, either initially or in the evolution, like for the solution to hyperbolic system of conservation laws. In this paper we identify the problem of lack of stability of the parareal in time algorithm and propose a simple way to cure it. The new method is used to solve a linear wave equation and a non linear Burger's equation, the results illustrate the stability of this variant of the parareal in time algorithm. \end{abstract}



\pagestyle{myheadings}
\thispagestyle{plain}
\markboth{X. Dai and Y. Maday}{ Stable parareal in time algorithm for   hyperbolic system}

\section{Introduction}
Parallel in time algorithms represent a new competitive way of using the ever
increasing number of cores in today's supercomputers platforms in cases where classical domain decomposition methods either are
of no relevance (for systems of differential equations) or in cases where the scalability of the sole domain decomposition (for evolution PDE's) gets to saturation. This kind of algorithm in the $4^{th}$ direction has received too little attention in the literature in the past due to the inherent sequential nature of the time direction: it seems indeed difficult to simulate what will arrive in the far future,
e.g. next week,  without knowing in details what will arrive in the close future, e.g. tomorrow and each day of the current week. Due to the paramount importance of finding new ways for parallelisation, going together with the ever increasing number of processors present in the new generations of architectures, algorithms such as
parareal in time, parallel deferred corrections or parallel exponential
integrators have become an active research topic
showing promising applicability for large-scale computations and perhaps one way to go to exascale computations. We refer e.g. to the book of K. Burrage~\cite{Bur:book} for
a synthetic approach on the subject (see also~\cite{burrage_review}). In
this book, the various
techniques of parallel in time algorithms are classified into three
categories: (i) parallelism {\em across the system}, (ii) parallelism
{\em across the method} and (iii) parallelism {\em across the time}.

The parareal in time algorithm pioneered in \cite{LMT}  extended later in \cite{BBMTZ;PTMDS, {BM;PTDNPDE}, MT;PTPCPDE} under a form much better tuned to the treatment of nonlinear problems  enters in the third category above (see \cite{Nievergelt} for a first algorithm in this category). It combines two propagators, one fine and one coarse, in a predictor-corrector manner that allows to use the fine propagator in parallel on the various processors, while the coarse propagator is solved  in a sequential manner. The algorithm provides good, provable convergence, for diffusive PDEs such as parabolic type problem. However, several studies have shown some instabilities for hyperbolic equations \cite{Farhat-2003, Fisher-2005,  mercerat:1521, chouly:1517} and the numerical analysis performed on simple examples or partial differential equations in \cite{{BM;PTDNPDE},Bal-2005,Gunnar} allows us to understand this behavior.

A challenging  field of application for parallel in time propagators is for long time simulation of Hamiltonian systems. Here also the plain method lacks the geometric properties  that are essential to guarantee the quality of the approximation. Let us remind that these geometric properties lead to the convergence of the algorithm towards the solution of the original systems on long time simulations. These properties ---either the symplecticity or the symmetry of the scheme --- are not preserved for the parareal  scheme even if the basics coarse and fine propagators are symplectic or symmetric.
A cure to this lack of geometric properties has been proposed in \cite{DLBLM}, based on two ingredients: the first one is the symmetrization of the parareal in time algorithm that leads to a new multistep schemes. This ingredient alone is not sufficient though to provide a correct algorithm for Hamiltonian systems. Indeed, some resonances are artificially introduced by the parareal strategy itself as is explained in \cite{DLBLM}; they prevent the symmetric variant of behaving well. The second ingredient  that has been introduced in \cite{DLBLM}
consists simply, after each correction of the predictor corrector scheme,  to project the solution coming out the parareal in time algorithm, on the manifold expressing the preservation of some basic invariant quantities of the Hamiltonian system. This is quite simple to implement, at least in its basic formulation (there exists a symmetric version based on a more involved implicit resolution of projection on the manifold).

Following this second idea advocated to stabilize the parareal in time algorithm for Hamiltonian system, we propose in this paper an extension of it  allowing to cure the lack of stability that occurs when solving hyperbolic type problem. Another strategy, based on providing a much more involved correction procedure re-utilizing previously computed information has been proposed in \cite{Farhat-2006, Farhat-2009} and already improved  the performance of the parareal in time algorithm for second-order hyperbolic systems. The method proposed here seems much more simple to implement and more general in scope and applications.

The remainder of this paper is organized as follows. In Section 2, we remind the basics on the plain parareal in time algorithm including the current state of the art on the stability and convergence analysis of the algorithm; we use this to introduce the proposed stable variant of the parareal scheme. In Section 3 we address the case of the wave equation, a simple second-order hyperbolic equation. The correction we propose is discussed in full details and the numerical simulation illustrates the good properties of the variant. In Section 4 we address the case of the Burgers' equation, with small viscosity, here again the variant performs well. In Section 5 we present some elements of numerical analysis proving  the stability of the proposed variant of the parareal scheme and finally some conclusions are drawn in Section 6.

\section{Some elements on the parareal in time algorithm}

\subsection{Basics on the parareal in time algorithm}

Let us consider a possibly non-linear  time dependent partial differential equation of the form
\begin{equation}
\label{eq:1}
\displaystyle\frac{\partial y}{\partial t} + {\mathcal A}(t;y) =
0, \qquad y(0) = y_0
\end{equation}
where $y(t)$ takes its values in a Banach space $B$ and $ {\mathcal A}$ is a
partial differential operator defined over $B$. We assume here that this problem is well posed, at least for some time:  the solution $y(t)$ exists for $t\in [0,T)$. More precisely, we assume that there exists
a propagator ${\mathcal E}$ such that
$$y(t) = {{\mathcal E}_{0}^{t}(y_{0})}, $$
named exact propagator. We note that, with obvious notations, we have
\begin{equation}
\label{eq2Kdix}
y(t) = {{\mathcal E}_{s}^{t}(y(s))},\quad \forall s, t, \quad 0\le s\le t\le T .
\end{equation}

Let us assume now that we want to simulate
this PDE on the interval $[0, T)$, we thus introduce a discrete propagator ${\mathcal F}$ that is assumed to be a fine and precise approximation of the exact solution based on a discretization in time only. The propagator
${\mathcal F}$  can be based on an approximation using  a time discretization and, e.g.  the use of an Euler-type or a more involved scheme based on a small enough time step $ \delta t$.

For instance, let us denote as $T_0=0$, $T_N=T$ and $T_n= n\Delta T$ (with $\Delta T = \frac{T}{N}$)
special times at which we are interested to consider snapshots of the
solution $y(T_n)$, then, with notations similar as above (\ref{eq2Kdix}), we have the approximation
$${Y_{n}} = {{\mathcal F}_{T_0}^{T_n}(y_{0})}, \quad n=1,...,N $$
we also have
$${Y_{n+1}} = {{\mathcal F}_{T_n}^{T_{n+1}}(Y_{n})}.$$

In what follows we present the parareal in time algorithm  that allows to define a sequence $\{{\bf Y}^k\}_k, {\bf Y}^k\equiv (Y_1^k, Y_2^k,...,Y_N^k)$, such that, for any $n$, ${\bf Y}_n^k $ converges to ${\bf Y}_n$  in the space $B$ as $k$ goes to
infinity. The definition of this sequence requires the availability of   another approximated propagator, denoted as  {${\mathcal G}$}. This second propagator is assumed to be cheaper than  ${\mathcal F}$ at the price of some
inaccuracy.
One can think immediately about having  {${\mathcal
G}$ based on a similar  Euler scheme as ${\mathcal F}$ with the larger time step $dT, \delta t<\!\!<dT\le \Delta T$}.
 But other possibilities are also offered: in addition to the solution scheme in time, a discretization in space has to be performed, for instance, the spacial discretization  can be based on a finite element method, the coarse solver ${\mathcal G}$ can then be based both on a coarser time step and a coarser spacial discretization, we refer to \cite{Fisher-2005} and \cite{MRS} where such a strategy is followed. Another possibility consists in providing
 ${\mathcal F}$
with all the physics of the phenomenon  while ${\mathcal G}$ will be based on a simplified physics (see eg \cite{Maday} for an example).

Assume ${\mathcal F}$ and ${\mathcal G}$ are given, we can  proceed to the definition of the sequence.
Let us initiate ${\bf Y}^0$ by
\begin{equation}\label{eq:00}
 Y^{0}_{n} =  {\mathcal G}_{T_0}^{T_n}(y_0), \quad n=1,...,N
\end{equation}
Assume now ${\bf Y}^k$ is known,  the iterative process proceeds by setting $Y^{k}_{0} =  y(0)$, and,
\begin{equation}Y^{k+1}_{n+1} = {\mathcal G}_{T_n}^{T_{n+1}}(Y^{k+1}_{n}) +{\mathcal F}_{T_n}^{T_{n+1}}(Y^{k}_{n}) -{\mathcal G}_{T_n}^{T_{n+1}}(Y^{k}_{n}), \quad n=1,...,N .
 \label{eq:predcorr}
\end{equation}

We first comment on the parallelization of the above algorithm. The initialization step (\ref{eq:00}) is sequential but involves the cheap coarse propagator. Then assuming that every $(Y^k_n)_n$ is known at step $k$, we first solve over each $]T_n,T_{n+1}[$ the totally independent propagation problems that provide  ${\mathcal F}_{T_n}^{T_{n+1}}(Y^{k}_{n})$; we can either propagate again (and now in parallel) the coarse solver to get ${\mathcal G}_{T_n}^{T_{n+1}}(Y^{k}_{n})$, or just remember that these quantities were computed during the previous step and stored, any way this is again parallel work. The corrector
${\mathcal F}_{T_n}^{T_{n+1}}(Y^{k}_{n}) - {\mathcal G}_{T_n}^{T_{n+1}}(Y^{k}_{n})$ can then be built, for any $n$ and the step can proceed. We finally  solve the predictor based on the coarse solver $ {\mathcal G}_{T_n}^{T_{n+1}}(Y^{k+1}_{n})$  --- sequentially ---   and immediately correct it by adding the previously computed corrections until convergence. Note that in this process, the fine solutions are only used on small propagation intervals and can be performed in parallel.

 \begin{figure}
   \begin{center}
\includegraphics[width = 8cm,angle =0]{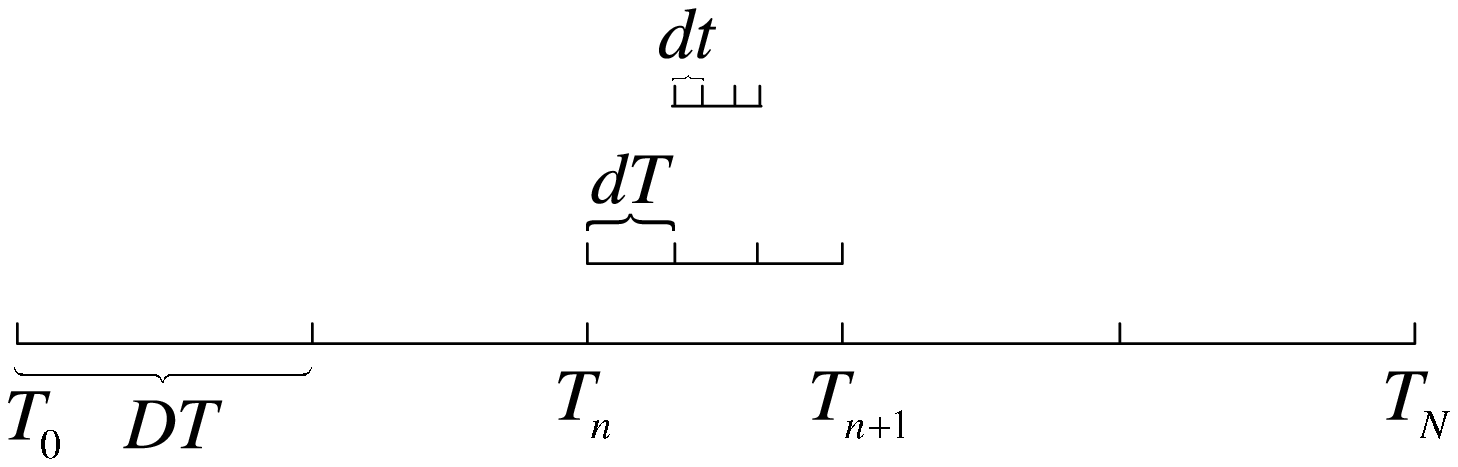}
\end{center}
\end{figure}

\subsection {Numerical analysis of the parareal in time algorithm}
We review in this subsection some elements on the stability and convergence of the parareal in time algorithm that have been first presented in \cite{Bal-2005}. Let us assume that there exists a sequence of Banach spaces $B_j$ such that $B_{j+1}\subset B_j$. It is also assumed that the problem \ref{eq:1} is stable in these spaces, in the sense that there exists a constant $C>0$, such that, for all $j, 1\le j\le J$
\begin{equation}
\label{eqno:02} \| {{\mathcal E}_{0}^{t}(y_{0})} \|_{B_j} \le C \| y_0\|_{B_j}
\end{equation}

The algorithm converges, under the following natural hypothesis: denoting by $\delta\mathcal{F}$ (resp. $\delta\mathcal{G}$) the difference $\delta\mathcal{F}_s^t = \mathcal{E}_s^t - \mathcal{F}_s^t $ (resp. $\delta\mathcal{G}_s^t = \mathcal{E}_s^t - \mathcal{G}_s^t $), we assume that there exists a constant $C>0$, such that, for all $j, 1\le j\le J$
\begin{equation}
\forall t \ge 0,  \forall \tau\ge 0,\quad
\|\mathcal{G}_t^{t+\tau}(x)  - \mathcal{G}_t^{t+\tau}(y)  \|_{B_j}\le  (1+ C \tau) ( \|x-y\|_{B_{j}}).\label{eq:1.4}
\end{equation}
In addition let us assume that there exists a constant $C>0$, such that, for all $j, 1\le j\le J$
\begin{eqnarray}
\forall t \ge 0,  \forall \tau\ge 0,\quad
\|\delta\mathcal{F}_t^{t+\tau}(x) - \delta\mathcal{F}_t^{t+\tau}(y)\|_{B_j}&\le& C \tau \eta\|x-y\|_{B_{j+1}},
 \label{eq:1.5}
\\   \|\delta\mathcal{G}_t^{t+\tau}(x) - \delta\mathcal{G}_t^{t+\tau}(y) \|_{B_j}&\le& C \tau \varepsilon \|x-y\|_{B_{j+1}}, \label{eq:1.6'}
\end{eqnarray}
where $\eta$ and $\varepsilon$ are small quantities, typically, for the Euler scheme these are $\eta=\delta t$ and $\varepsilon = \Delta T$ and $C$ is a generic constant, independent of the quantities on the right hand side of the previous inequalities and on $j\le J$, but may depend on $J$.

\begin{remark} For Euler schemes based on a time step $\delta t$ for $\mathcal{F}$ and $\Delta T$ $\mathcal{G}$, hypothesis (\ref {eq:1.5}) and (\ref {eq:1.6'}) hold with $\eta=\delta t$ and $\varepsilon = \Delta T$ respectively.
\end{remark}

The convergence is stated in the following (see theorem 1 in  \cite{Bal-2005})
\begin{theorem}\label{theo:error}
  Assume that the discrete propagators $\mathcal{F}$ and $\mathcal{G}$ satisfy (\ref{eq:1.4}), (\ref{eq:1.5}) and (\ref{eq:1.6'}), then the error between the exact solution and the solution provided by the parareal scheme (\ref{eq:predcorr}) satisfies for any $k\ge 0$:
 \begin{equation}
\forall n\le N,\quad \| Y_{n}^{k} - y(T^n)\|_{B_0} \le C(\varepsilon^{k+1} + \eta)(1+ \|y_0\|_{B_k}).
 \label{eq:1.6}
\end{equation}
  \end{theorem}

\noindent {\bf Proof}: The proof is done by induction over $j\ge 0$ for the statement

 \begin{equation}
\forall i\ge 0, \quad \| Y_{n}^{j} - y(T^n)\|_{B_i} \le C(\varepsilon^{j+1} + \eta)(1+ \|y_0\|_{B_{i+j+1}}),
 \label{eq:1.7}
\end{equation}
which is obvious for $j=0$ since it corresponds to the analysis of the plain coarse operator. Let us assume that it is true for $j$ and let us prove it for $j+1$.
From the definition of $Y_{n+1}^{j+1}$, we write
\begin{eqnarray*}
Y_{n+1}^{j+1} - y(T^{n+1}) &=&  {\mathcal G}_{T_n}^{T_{n+1}}(Y^{j+1}_{n}) +{\mathcal F}_{T_n}^{T_{n+1}}(Y^{j}_{n}) -{\mathcal G}_{T_n}^{T_{n+1}}(Y^{j}_{n}) - {\mathcal E}_{T_n}^{T_{n+1}}[y(T^n) ]\\
&=& {\mathcal G}_{T_n}^{T_{n+1}}(Y^{j+1}_{n}) - {\mathcal G}_{T_n}^{T_{n+1}}(y(T^n)) \\
&& + [{\mathcal E}_{T_n}^{T_{n+1}} - {\mathcal G}_{T_n}^{T_{n+1}}] (Y^{j}_{n}) - [{\mathcal E}_{T_n}^{T_{n+1}} - {\mathcal G}_{T_n}^{T_{n+1}}] (y(T^n))\\
&& \quad- [{\mathcal E}_{T_n}^{T_{n+1}} - {\mathcal F}_{T_n}^{T_{n+1}}] (Y^{j}_{n}) + [{\mathcal E}_{T_n}^{T_{n+1}} - {\mathcal F}_{T_n}^{T_{n+1}}] (y(T^n))\\
&& \qquad - [{\mathcal E}_{T_n}^{T_{n+1}} - {\mathcal F}_{T_n}^{T_{n+1}}] (y(T^n))
\end{eqnarray*}
or again with the notations introduced above
\begin{eqnarray*}
Y_{n+1}^{j+1} - y(T^{n+1}) &=& {\mathcal G}_{T_n}^{T_{n+1}}(Y^{j+1}_{n}) - {\mathcal G}_{T_n}^{T_{n+1}}(y(T^n)) \\
&& +\delta {\mathcal G}_{T_n}^{T_{n+1}} (Y^{j}_{n}) - \delta {\mathcal G}_{T_n}^{T_{n+1}} (y(T^n))\\
&& \quad- \delta {\mathcal F}_{T_n}^{T_{n+1}} (Y^{j}_{n}) + \delta  {\mathcal F}_{T_n}^{T_{n+1}} (y(T^n))\\
&& \qquad - [{\mathcal E}_{T_n}^{T_{n+1}} - {\mathcal F}_{T_n}^{T_{n+1}}] (y(T^n))
\end{eqnarray*}
hence by triangular inequality
\begin{eqnarray*}
 \| Y_{n+1}^{j+1} - y(T^{n+1})\|_{B_i} &\le& \|{\mathcal G}_{T_n}^{T_{n+1}}(Y^{j+1}_{n}) - {\mathcal G}_{T_n}^{T_{n+1}}(y(T^n))\|_{B_i} \\
&& +\|\delta {\mathcal G}_{T_n}^{T_{n+1}} (Y^{j}_{n}) - \delta {\mathcal G}_{T_n}^{T_{n+1}} (y(T^n))\|_{B_i}\\
&& \quad+\| \delta {\mathcal F}_{T_n}^{T_{n+1}} (Y^{j}_{n}) + \delta  {\mathcal F}_{T_n}^{T_{n+1}} (y(T^n))\|_{B_i}\\
&& \qquad +\| [{\mathcal E}_{T_n}^{T_{n+1}} - {\mathcal F}_{T_n}^{T_{n+1}}] (y(T^n))\|_{B_i}
\end{eqnarray*}
that from hypothesis (\ref{eq:1.4}), (\ref{eq:1.5}) and (\ref{eq:1.6'}) can be upper bounded as follows
\begin{eqnarray*}
 \| Y_{n+1}^{j+1} - y(T^{n+1})\|_{B_i} &\le&(1+C\Delta T) \| Y^{j+1}_{n} - y(T^n)\|_{B_i}\\
&& + C \Delta T \varepsilon  \| Y^{j}_{n} - y(T^n)\|_{B_{i+1}}\\
&& \quad+ C \Delta T \eta \| Y^{j}_{n} - y(T^n)\|_{B_{i+1}}   + C \Delta T \eta \|y(T^n)\|_{B_{i+1}}
\end{eqnarray*}
the stability hypothesis (\ref{eqno:02}) together with the induction hypothesis (\ref{eq:1.7}) allow to get
\begin{eqnarray*}
 \| Y_{n+1}^{j+1} - y(T^{n+1})\|_{B_i} &\le&(1+C\Delta T) \| Y^{j+1}_{n} - y(T^n)\|_{B_i}\\
&& + C \Delta T (\varepsilon  +  \eta) (\varepsilon^{j+1} +\eta)  \| y_0\|_{B_{i+j+2}}   + C \Delta T \eta \|y_0\|_{B_{i+1}}\\
&\le&(1+C\Delta T) \| Y^{j+1}_{n} - y(T^n)\|_{B_i}\\
&& \quad + C \Delta T  (\varepsilon^{j+2} + \eta) \| y_0\|_{B_{i+j+2}}
\end{eqnarray*}
from which we get  (\ref{eq:1.7}) for $j+1$ from the discrete Gronwall lemma, and the proof is complete.

\bigskip
It is interesting to notice two things at this level. First the constant $C$ in (\ref{eq:1.6}) depends on $k$ and $T$ and this convergence result can be polluted on a long time range by a constant $C$ that becomes too large. Second, the convergence results implies, of course, stability. Hence, provided that  the initial condition is regular enough and remains regular the parareal in time algorithm is stable.

The previous remark explains actually the problem of the use of the parareal scheme when non dissipative PDE or PDE with small dissipation are simulated. This is the case for small viscosity parabolic problem as is explained in the next subsection.

This is also the case for hyperbolic systems. We can illustrate this statement on two equations: the periodic wave equation and the periodic Burger's equation for which, starting with regular enough initial condition solves or corrects the instability of the parareal in time algorithm in this situation. This allows to rectify somehow the statements that were done before in e.g. \cite{Farhat-2009} where the stability problem occurs for second-order hyperbolic problems and not for  parabolic and first-order
hyperbolic problems. The problem is not really on the type of the equations (and certainly not on the order of the equations) but on the regularity of the solution, including the regularity of the initial condition.

\subsection{Behavior on a parabolic case}\label{parab}

We have performed a series of simple test on the linear parabolic equation, provided with periodic boundary conditions
\begin{eqnarray}\label{parabolic-1d}
\left\{\begin{array}{rcl}
\frac{\partial u}{\partial t} - \mu \frac{\partial^2 u}{\partial x^2} +  \frac{\partial u}{\partial x}  &=& 0, \\
u(x,0) &=& u_0(x),  \\
u(x,t) &=& u(x+2\pi,t),
\end{array}
\right.
\end{eqnarray}
in the case where the viscosity $\mu$ is positive yet very small $\mu = 10^{-10}$. We have chosen a standard spectral Fourier approximation in space (truncated series with $N=256$ and chosen a second order in time implicit Euler scheme
for a propagation over $[0,T=2]$. The parareal in time algorithm is implemented with $\Delta T= 2. 10^{-2}$, $dT= 10^{-2}$ and $\delta t = 10^{-4}$.

The simulation are done in order to illustrate the behavior of the scheme as regards the regularity of the solution. We have thus chosen different types of initial conditions with different regularities by choosing $u_0(x) = \sum \hat u_\ell e^{i\ell x}$ with

\begin{eqnarray*}
\hat u_\ell = \left\{\begin{array}{rcl}
&\frac{0.5 i}{|\ell|^p}, \quad &\hbox{if  $ \ell>0 $}  \\
&0, \quad  & \hbox{if $\ell=0$}\\
&- \frac{0.5 i}{|\ell|^p}, \quad  &\hbox{if  $\ell<0$},
\end{array}
\right.
\end{eqnarray*}

In the case of low regularity: $p=0.5$ or $p=1$, the relative error after 15 parareal iterations is about 100! There is no blow up, the iterations can continue, the error eventually going to zero. For $p=4$ or even better
$p=10$, the convergence is achieved after very few iterations, as is expected from Theorem \ref{theo:error}.

These simple computations provide a preliminary illustration that the behavior of the parareal in time algorithm is deeply linked to the regularity of the solution.

\section{The periodic wave equation}

\subsection{The basic Fourier discretization}

We consider the following one dimensional  wave equation with periodic boundary conditions: Find $u$ such that
\begin{eqnarray}\label{wave-eqaution-1d}
\left\{\begin{array}{rcl}
\frac{\partial^2 u}{\partial t^2} - c^2 \frac{\partial^2 u}{\partial x^2}&=& 0, \\
u(x,0) &=& f(x),~~  \frac{\partial u}{\partial t}(x, 0)  = g(x), \\
u(x,t) &=& u(x+\mathbb{T},t),~~ \frac{\partial u}{\partial x}(x,t) =  \frac{\partial u}{\partial x}(x+\mathbb{T},t)
\end{array}
\right.
\end{eqnarray}
where $f$ and $g$ are periodic given initial condition, $c$ stands for the speed of the wave, and $\mathbb{T}\in \RR^+$ represents the period.

The weak form of equation (\ref{wave-eqaution-1d}) is as follows: Find $u$ such that $\forall t>0, u(.,t) \in S = \{ v\in H^1(0,\mathbb{T}), v(0)=v(\mathbb{T})\}$
\begin{eqnarray} \label{weak-form-wave-eqaution-1d}
(\frac{\partial^2 u}{\partial t^2}, v) + c^2(\frac{\partial u}{\partial x}, \frac{\partial v}{\partial x}) = 0, \quad \forall v
\in  S,
\end{eqnarray}
here,  $(u, v) = \int_{0}^{\mathbb{T}} u \bar{v} dx$.

We use now the Fourier spectral method to semi discretize in space
problem (\ref{weak-form-wave-eqaution-1d}). That is, we choose
\begin{eqnarray}\label{space-spectral}
S_N = span\{ e^{i\ell \frac{2\pi}{\mathbb{T}}x}, \ell=-N, \cdots,
N \},
\end{eqnarray}
 and approximate $u(x, t)$ by
\begin{eqnarray}\label{approx-u}
u_N(x, t) = \sum_{\ell=-N}^{N} \hat{u}_{\ell}(t)
e^{i\ell\frac{2\pi}{\mathbb{T}}x}.
\end{eqnarray}
Note that $\hat{u}_\ell(t) = \frac{1}{\mathbb{T}}
(u_N(x, t), e^{i\ell\frac{2\pi}{\mathbb{T}}x})$, as results from the equality  $(e^{i\ell\frac{2\pi}{\mathbb{T}}x},
e^{ik\frac{2\pi}{\mathbb{T}}x}) = \mathbb{T} \delta_{\ell,k}$ where $ \delta_{\ell,k}$ denotes the standard Kroneker symbol.

The solution $u_N$ is defined by a Galerkin approximation so as to satisfy
\begin{eqnarray} \label{weak-form-wave-eqaution-1d_N}
(\frac{\partial^2 u_N}{\partial t^2}, v_N) + c^2(\frac{\partial u_N}{\partial x}, \frac{\partial v_N}{\partial x}) = 0, \quad \forall v_N
\in  S_N.
\end{eqnarray}
As is classical, this Galerkin approximation reduces to a system of differential equations in the Fourier coefficients $\hat{u}_{k}(t)$. More precisely, let us choose iteratively $v_N=e^{ik\frac{2\pi}{\mathbb{T}}x}$, $k \in [-N, N]$ in (\ref{weak-form-wave-eqaution-1d_N}).
We then obtain
the following system of independent equations
\begin{eqnarray}
&&\mathbb{T} \frac{\partial^2 \hat{u}_{\ell}}{\partial t^2}(t) +
\frac{(2\pi)^2}{\mathbb{T}} \ell^2 c^2 \hat{u}_{\ell}(t) = 0,
~ \forall \ell \in [-N, N ],\label{mode_ell}
\end{eqnarray}
that can be exactely solved.
After resolution, we know all the coefficients  $\hat{u}_\ell(t)$ in a closed form,
then from equation   (\ref{approx-u})
we get the Galerkin approximation of $u(x, t)$ in $S_N$.

Next, using respectively $v = \frac{\partial u}{\partial t}$ in (\ref{weak-form-wave-eqaution-1d}) and  $v_N = \frac{\partial u_N}{\partial t}$ in (\ref{weak-form-wave-eqaution-1d_N})
allows to derive, for any time $t$
\begin{equation}
\| \frac{\partial u}{\partial t}(.,t)\|^2_{L^2(0,\mathbb{T})} + c^2 \| \frac{\partial u}{\partial x}(.,t)\|^2_{L^2(0,\mathbb{T})} = \| \frac{\partial u}{\partial t}(.,0)\|^2_{L^2(0,\mathbb{T})} + c^2 \| \frac{\partial u}{\partial x}(.,0)\|^2_{L^2(0,\mathbb{T})}
\end{equation}
and
\begin{equation}\label{conservation_wave}
\| \frac{\partial u_N}{\partial t}(.,t)\|^2_{L^2(0,\mathbb{T})} + c^2 \| \frac{\partial u_N}{\partial x}(.,t)\|^2_{L^2(0,\mathbb{T})} = \| \frac{\partial u_N}{\partial t}(.,0)\|^2_{L^2(0,\mathbb{T})} + c^2 \| \frac{\partial u_N}{\partial x}(.,0)\|^2_{L^2(0,\mathbb{T})}
\end{equation}
respectively. If we now introduce the following Hamiltonian
\begin{equation}\label{definition_H}
H(u_N)(t) = \| \frac{\partial u_N}{\partial t}(.,t)\|^2_{L^2(0,\mathbb{T})} + c^2 \| \frac{\partial u_N}{\partial x}(.,t)\|^2_{L^2(0,\mathbb{T})}
\end{equation}
The conservation (\ref{conservation_wave}) reads for any time $t$
\begin{equation}
H(u_N)(t) = H(u_N)(0)
\end{equation}

Of course the numerical discretization in time should preserve this invariant .

\begin{remark}\label{remark1}
As is often the case for Hamiltonian system, this basic invariant (which is here the total energy for the wave equation) is not the only one available, for instance it is immediately seen on (\ref{mode_ell}) that  by defining
\begin{equation}
h_\ell(v) = | \frac{\partial v}{\partial t} |^2 + (\frac{2\pi}{\mathbb{T}})^2  \ell^2 c^2 | v |^2,
\end{equation}
for any $\ell$, $\ell=-N,..,N$, we have
\begin{equation}
\forall t, \quad h_\ell(\hat u_\ell)(t) = h_\ell(\hat u_\ell)(0)
\end{equation}
which provides $2N+1$ invariant quantities (we notice that $H= \sum_{\ell = -N}^N h_\ell$).
\end{remark}

\subsection{Illustration of the stability and instabilities of the plain parareal in time algorithm}

In all what follows, the discrete solutions will always, for any time $t$, belong to $S_N$ and will result from a  Galerkin process. We will thus denote by $u_n$ the approximation (in time) of the value of $u_N(T_n)$ (hence the full approximation of $u(T_n)$).
In order to test the parareal in time algorithm on this simple second order equation, we further discretize the semi discrete problem (\ref{weak-form-wave-eqaution-1d_N}) by adding a time discretization, that is based, both for the coarse $ \mathcal{G}_{\Delta T}$ and the fine $ \mathcal{F}_{\Delta T}$ propagators, on the velocity Verlet scheme which is a second order in time, stable and consistent with the preservation of invariants (the wave equation being an Hamiltonian system, the symplectic nature of the velocity Verlet scheme makes it a method of choice).

The plain parareal  scheme reads:
\begin{eqnarray}
\left\{\begin{array}{rcl}
  u_{n+1}^{0} &=& \mathcal{G}_{\Delta T}(u_n^0) \\
  u_{n+1}^{k+1} &=& \mathcal{G}_{\Delta T}(u_{n}^{k+1}) +  \mathcal{F}_{\Delta T}(u_{n}^{k})
                  - \mathcal{G}_{\Delta T}(u_{n}^{k}).\label{para1}
  \end{array}
\right.
 \end{eqnarray}

We perform simulations with the discretization parameters equal to
\begin{enumerate}
\item $T = 100, \Delta T = 2. 10^{-1}, dT = 4. 10^{-3}, \delta t = 2. 10^{-6}$.
\item N = 30.
\end{enumerate}
for different  initial conditions of different regularities: $u_0(x) = \sum_{\ell\in Z} \hat u _\ell  e^{i\ell\frac{2\pi}{\mathbb{T}}x}$
with

\begin{eqnarray*}
\hat u_\ell = \left\{\begin{array}{rcl}
&\frac{1 }{|\ell|^p}, \quad &\hbox{if  $ \ell\not= 0 $}  \\
&0, \quad  & \hbox{if $\ell=0$},
\end{array}
\right.
\end{eqnarray*}

with p
$p=1, 5$ or $10$. The cases with p = 1 and p = 5 diverges while the case with
p = 10 converges.

\subsection{Results with the plain parareal in time algorithm on a realistic test problem}

In this example, the excitation is provided by a Ricker pulse, which is analytically defined as follows

 \begin{enumerate}
\item $\mathbb{T} = 5000$ m, $c = 2000$ m/s
\item and the initial conditions are $f(x) = 0$,  $g(x) = (1 - 2 (f_s \pi \frac{x-x_s}{c})^2) \exp(-(f_s \pi
\frac{x-x_s}{c})^2)$.
\end{enumerate}
Here, $f_s= 2.5$ Hz, is the frequency  and $x_s=2500$ m is the position of the vibration source.

 The numerical results are obtained with a coarse propagator $\mathcal{G}$ and
fine propagator $\mathcal{F}$ that differ only with time step $dT$ and $\delta t$, respectively :
\begin{enumerate}
\item $T = 100, \Delta T = 2. 10^{-1}, dT = 4. 10^{-3}, \delta t = 2. 10^{-6}$.
\item N = 30.
\end{enumerate}

\begin{figure}[!htbp]
\begin{center}
\includegraphics[width = 8cm, angle=-90]
{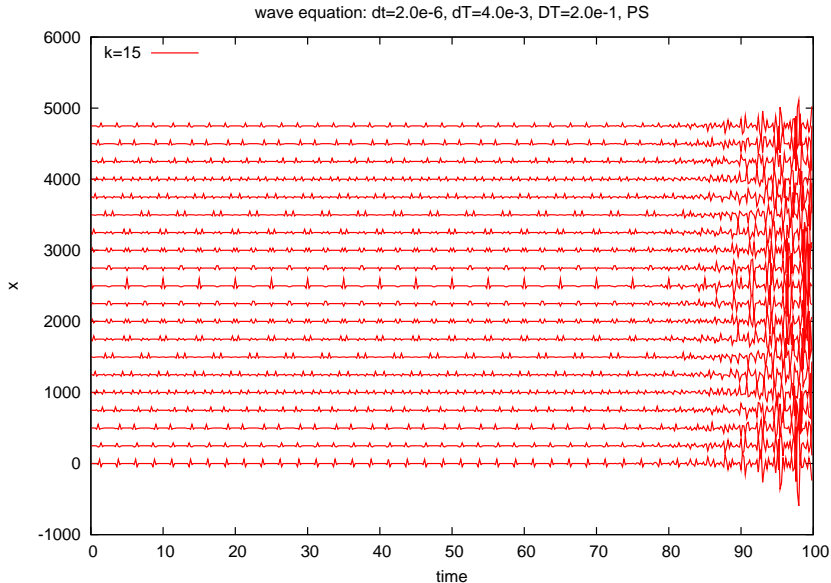}
\caption{Vibration  at different positions obtained by the parareal
method($\delta t = 2. 10^{-6}, dT=4. 10^{-3}, \Delta T=2.
10^{-1}$) revealing the instabilities.}\label{fig-wave-case1-2-PS-vibr-15}
\end{center}
\end{figure}

The plot showing the evolution in time of the solution obtained after 15 iterations at different positions in the periodic space $[0,,\mathbb{T}]$ on figure \ref{fig-wave-case1-2-PS-vibr-15} reveals an instability visible starting at time $t=80$. The situation is even worse after 25 iterations. Actually this is one case that illustrates the bad behavior of the parareal in time algorithm for this hyperbolic equation. Note however that for many other choices of the discretization parameters $\delta t$ or $dT$ or even $N$ but the very same initial condition, the simulation is perfectly fine, stable and convergent. The instability is thus not systematic, and changing slightly the discretization parameters, is often enough to make the scheme work. The plain parareal in time algorithm does lack of robustness, and this is what we want to correct in what follows.

\begin{figure}[!htbp]
\begin{center}
\includegraphics[width = 8cm, angle=-90]
{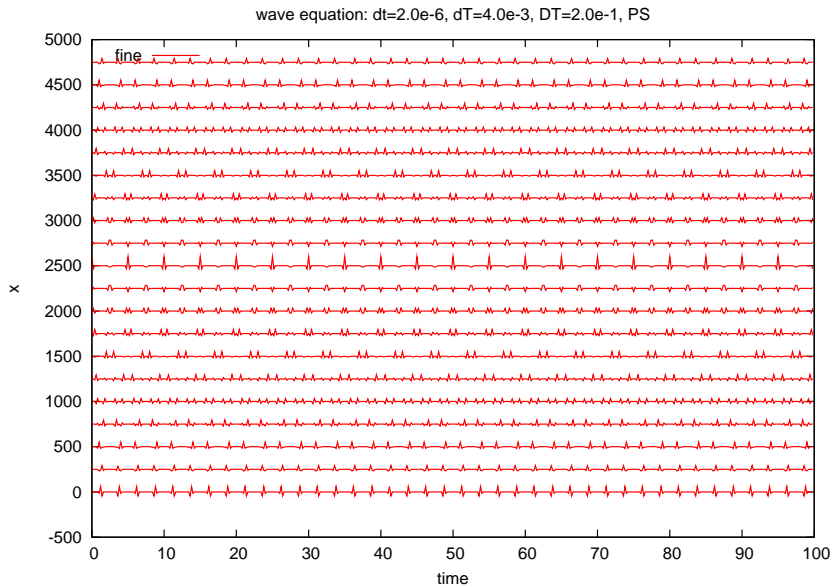}
\caption{Vibration  at different positions obtained by the fine
sequential
 method($\delta t = 2. 10^{-6}$)}\label{fig-wave-case1-2-vibr-fine}
\end{center}
\end{figure}

The exact solution shown on the next plot (figure
\ref{fig-wave-case1-2-vibr-fine}) illustrates the differences.

\begin{figure}[!htbp]
\begin{center}
\includegraphics[width = 8cm, angle=-90]
{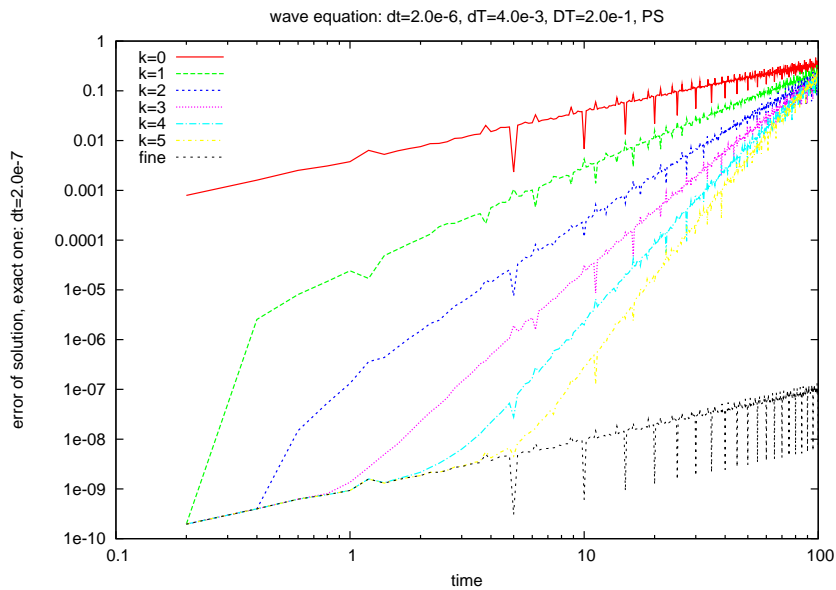}
\caption {Error on the solution obtained by the parareal
method($\delta t = 2. 10^{-6}, dT=4. 10^{-3}, \Delta T=2. 10^{-1}
$)}\label{fig-wave-case1-2-PS-P-0-5}
\end{center}
\end{figure}

\begin{figure}[!htbp]
\begin{center}
\includegraphics[width = 8cm, angle=-90]
{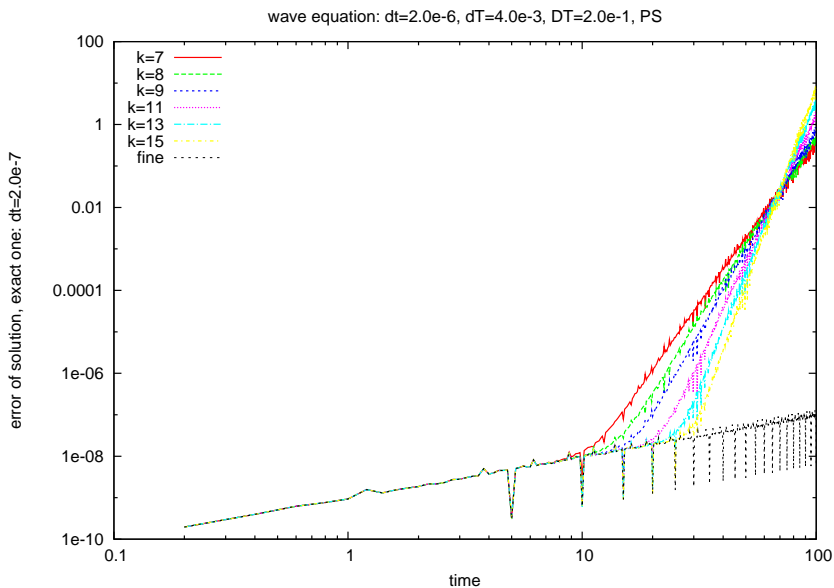}
\caption{Error on the solution obtained by the parareal
method($\delta t = 2. 10^{-6}, dT=4. 10^{-3}, \Delta T=2. 10^{-1}
$)}\label{fig-wave-case1-2-PS-P-7-15}
\end{center}
\end{figure}

The plots of the error between the exact solution and the solution obtained after $k$ iterations, with $k=1,..,15$ (see figures \ref{fig-wave-case1-2-PS-P-0-5} and \ref{fig-wave-case1-2-PS-P-7-15}) illustrate clearly that the solution obtained by the parareal in time algorithm
\begin{enumerate}
\item converges for short time, e.g. only 7 iterations are sufficient to recover the solution obtained with the sole sequential fine propagator at time equal to 10.
\item does not converge for longer propagations, here for a time equal to 100
\item During the first 5 to 6 iterations the error of the parareal in time algorithm is dominated by the error of the coarse scheme, but after $k=6$, the error increases largely to values that are order 1 and bigger
\item After 15 iterations of the parareal in time algorithm the parareal solution is computed with an acceptable error until time $t\simeq 30$  only.
\end{enumerate}

Taking into account the variants of the parareal in time algorithm introduced in \cite{DLBLM} to simulate  Hamiltonian systems  we infer that the bad behavior of the parareal scheme comes from the non preservation of the invariant quantities. Let us for example observe the evolution of the  value of the quantity $H$ defined in (\ref{definition_H})

\begin{figure}[!htbp]
\begin{center}
\includegraphics[width = 8cm, angle=-90]
{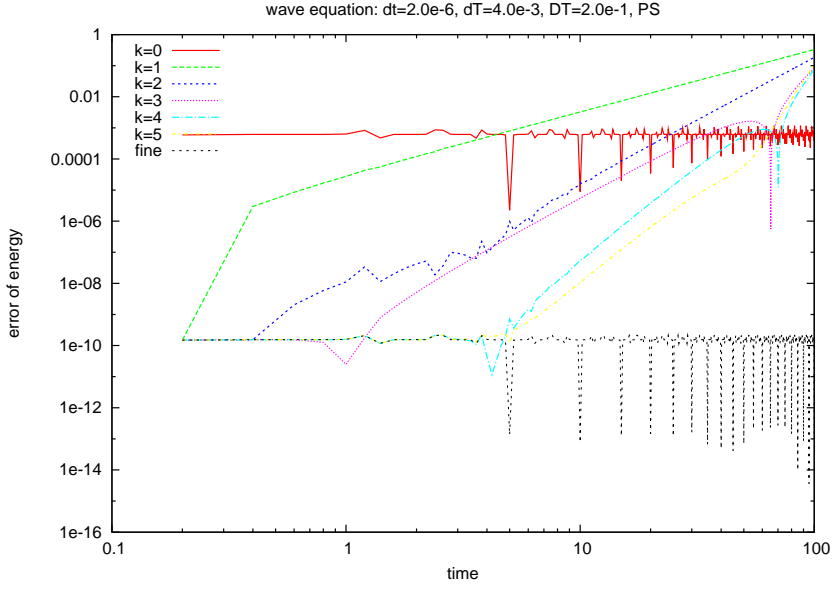}
\caption{Error on the energy obtained by the parareal method($\delta
t = 2. 10^{-6}, dT=4. 10^{-3}, \Delta T=2. 10^{-1}
$)}\label{fig-wave-case1-2-PS-E-0-5}
\end{center}
\end{figure}

\begin{figure}[!htbp]
\begin{center}
\includegraphics[width = 8cm, angle=-90]
{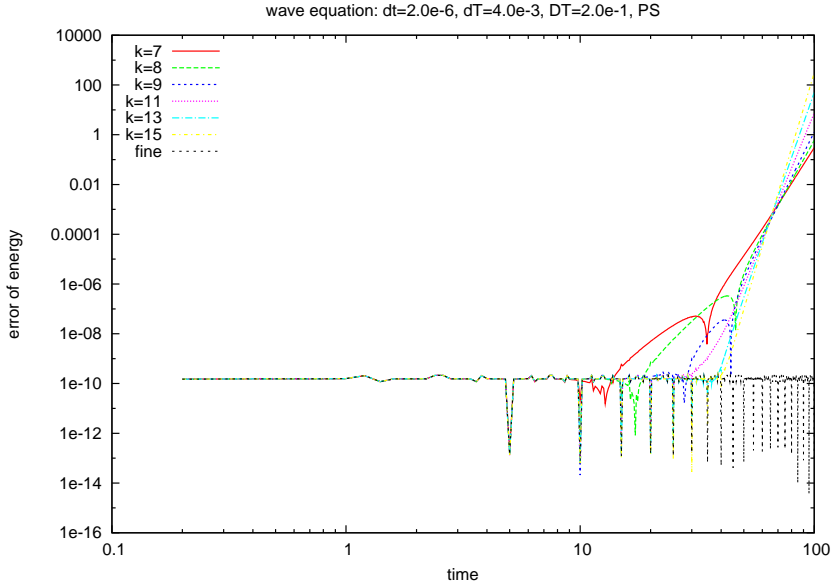}
\caption{Error on the energy obtained by the parareal method($\delta
t = 2. 10^{-6}, dT=4. 10^{-3}, \Delta T=2. 10^{-1}
$)}\label{fig-wave-case1-2-PS-E-7-15}
\end{center}
\end{figure}

The figures \ref{fig-wave-case1-2-PS-E-0-5} and \ref{fig-wave-case1-2-PS-E-7-15} illustrate the lack of conservation, that is coherent with   the divergence of the scheme in this case.

Finally we want to notice that we have performed the same simulation with a Ricker pulse of varying force (by modifying the value of $f_s$). The instability starts much earlier for larger $f_s$ confirming the relation between the instability and the singularity of the solution.

\subsection{Parareal algorithm with projection}

Let us define the manifold ${\cal M}$ corresponding to the set of all functions $v$ such that $H(v)=H(u(0))$.
The Standard Projection Method on ${\cal M}$  (see chapter 4 of
\cite{hairer-lubich-wanner}) can be written as
follows: For a given $\tilde u$, solve the system in $(u,\lambda)$
\begin{eqnarray*}
 u =   \widetilde{u} + \lambda
\nabla H(\widetilde{u}),
\end{eqnarray*}
where $\lambda \in
\mathbb{R}$  is such that  $H(u) = H(u(0))$, i.e.
\begin{eqnarray}\label{eq:proj-eq-1}
H(\widetilde{u} + \lambda \nabla H(\widetilde{u})) = H(u(0)).
\end{eqnarray}
We then set
\begin{equation}
u = \pi_{\cal M}(\widetilde{u})
\end{equation}

In the general case, equation (\ref{eq:proj-eq-1}) can only be solved  by some iterative
method, for example, Newton iterations, but in the current case, using the fact that $H$ is quadratic, we notice and take into account that equation (\ref{eq:proj-eq-1}) is a second order polynomial
in $\lambda$ that can be solved exactly.

Then, combining the  plain parareal  scheme (\ref{para1})
with the projection method above,  we obtain the
parareal method  with projection to ${\cal M}$:
\begin{eqnarray} \label{proj_para}\left\{
\begin{array}{rcl}
  u_{n+1}^{0} &=& \mathcal{G}_{\Delta T}(u_n^0) \\
\tilde{u}_{n+1}^{k+1} &=& \mathcal{G}_{\Delta T} (u_{n}^{k+1}) +
\mathcal{F}_{\Delta T} (u_{n}^{k}) - \mathcal{G}_{\Delta T}
(u_{n}^{k}),
\\
u_{n+1}^{k+1} &=& \pi_{\cal M}(\tilde{u}_{n+1}^{k+1}).
\end{array}
\right.
\end{eqnarray}
  
\begin{figure}[!htbp]
\begin{center}
\includegraphics[width = 8cm, angle=-90]
{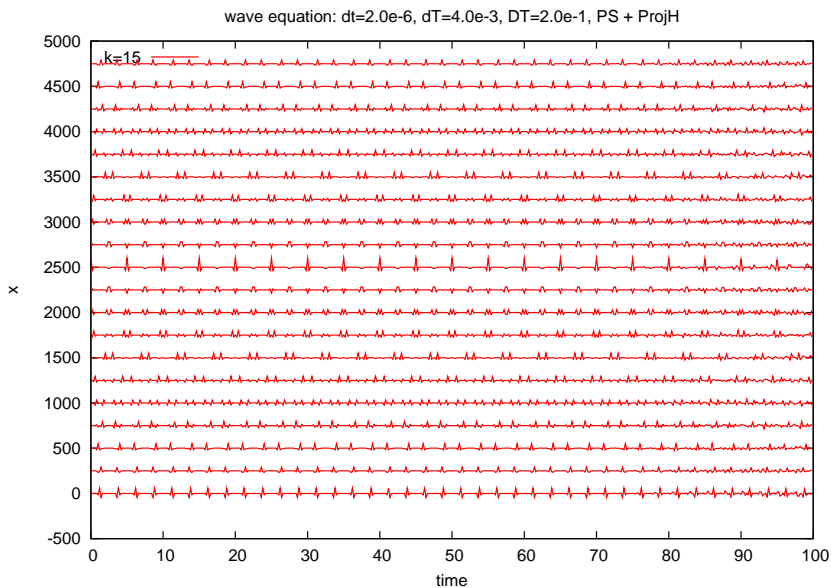}
\caption{Vibration  at different positions obtained by the parareal
method with projection to the manifold ${\cal M}$ ($\delta t = 2.
10^{-6}, dT=4. 10^{-3}, \Delta T=2. 10^{-1}
$)}\label{fig-wave-case1-2-PS-ProjH-vibr-15}
\end{center}
\end{figure}

 We now show the results
obtained by parareal method with projection. We first note on figure
\ref{fig-wave-case1-2-PS-ProjH-vibr-15} that the global solution
after $k=15$ iterations does not seem to present the instabilities
of the plain parareal  in time algorithm. The same conclusion holds after 25 iterations.

\begin{figure}[!htbp]
\begin{center}
\includegraphics[width = 8cm, angle=-90]
{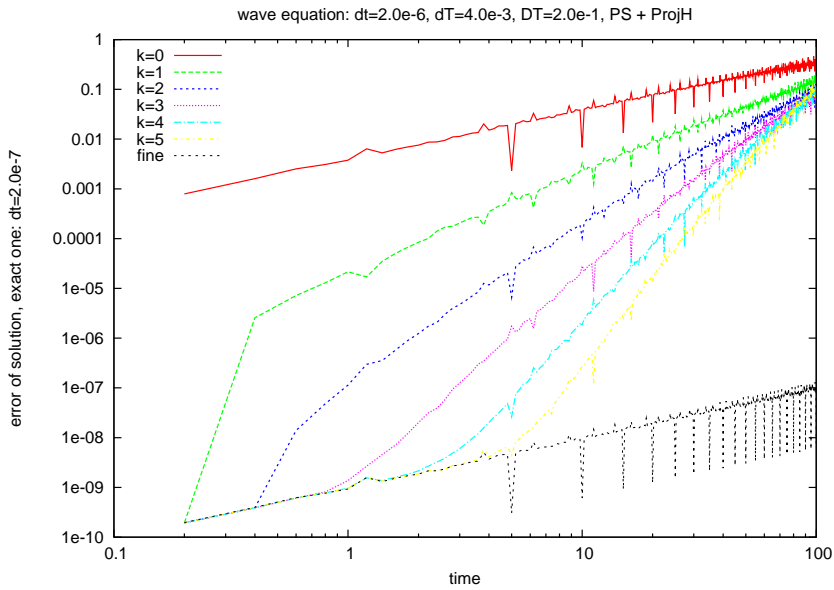}
\caption{Error on the solution obtained by the parareal method with
projection to the manifold ${\cal M}$ ($\delta t = 2. 10^{-6},
dT=4. 10^{-3}, \Delta T=2. 10^{-1}
$)}\label{fig-wave-case1-2-PS-ProjH-P-0-5}
\end{center}
\end{figure}

\begin{figure}[!htbp]
\begin{center}
\includegraphics[width = 8cm, angle=-90]
{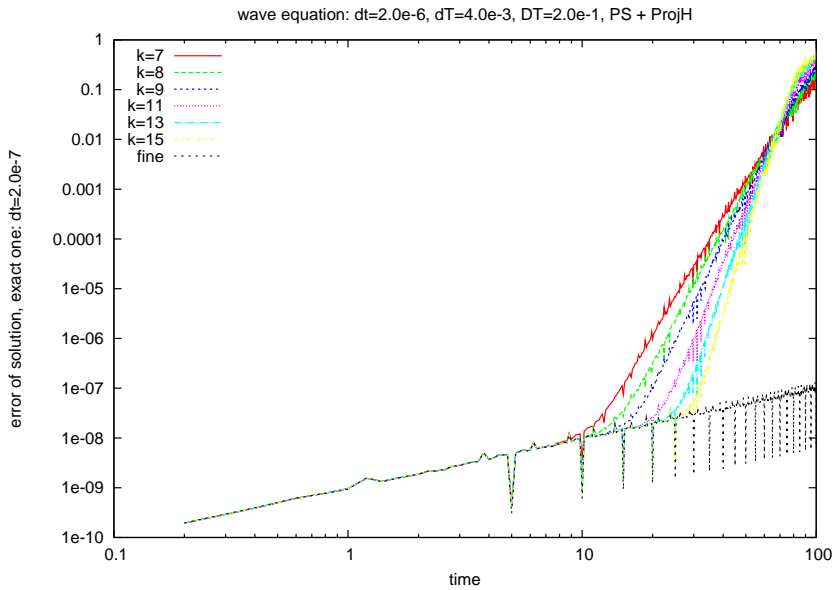}
\caption{Error on the solution obtained by parareal method with
projection to the manifold ${\cal M}$ ($\delta t = 2. 10^{-6},
dT=4. 10^{-3}, \Delta T=2. 10^{-1}
$)}\label{fig-wave-case1-2-PS-ProjH-P-7-15}
\end{center}
\end{figure}

This is confirmed by the figures \ref{fig-wave-case1-2-PS-ProjH-P-0-5} and \ref{fig-wave-case1-2-PS-ProjH-P-7-15} where we note an improvement in the error between the exact solution and the discrete one as the iterations evolve. The error is always at worse of the same order as the error of the coarse propagator.

Of course this is far from being what is expected, but the conservation of $H$ does improve the behavior of the algorithm. The following figures
\ref{fig-wave-case1-2-PS-ProjH-E-0-5}
and
\ref{fig-wave-case1-2-PS-ProjH-E-7-15} illustrate the good preservation of $H$ along the trajectory for the parareal method with projection. Although the solution does not diverge as in the  plain
parareal case, it does not converge as we would like.

\begin{figure}[!htbp]
\begin{center}
\includegraphics[width = 8cm, angle=-90]
{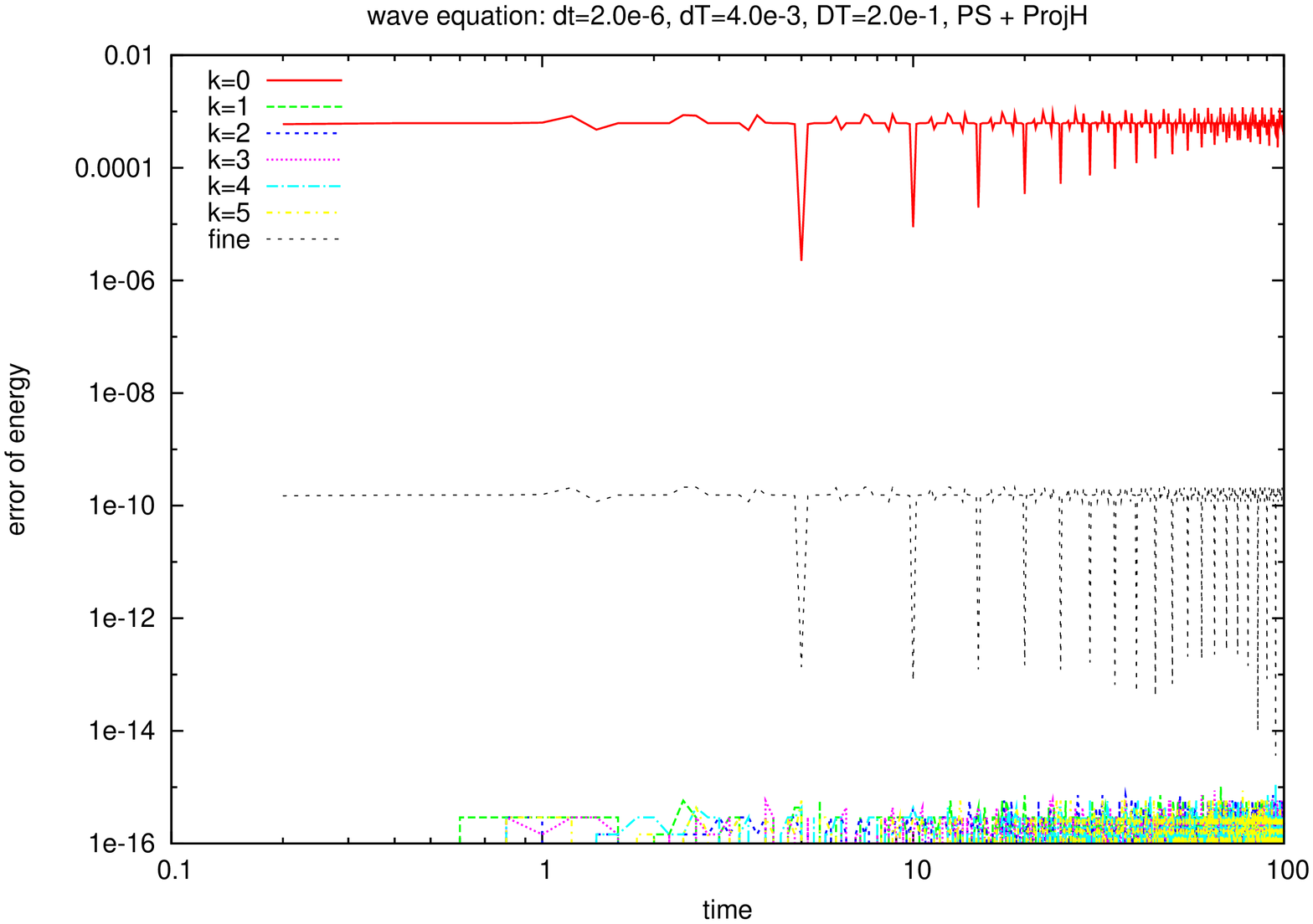}
\caption{Error on the energy obtained by
 parareal method with projection to the manifold ${\cal M}$ ($\delta t = 2. 10^{-6}, dT=4. 10^{-3}, \Delta T=2. 10^{-1} $)}\label{fig-wave-case1-2-PS-ProjH-E-0-5}
\end{center}
\end{figure}

\begin{figure}[!htbp]
\begin{center}
\includegraphics[width = 8cm, angle=-90]
{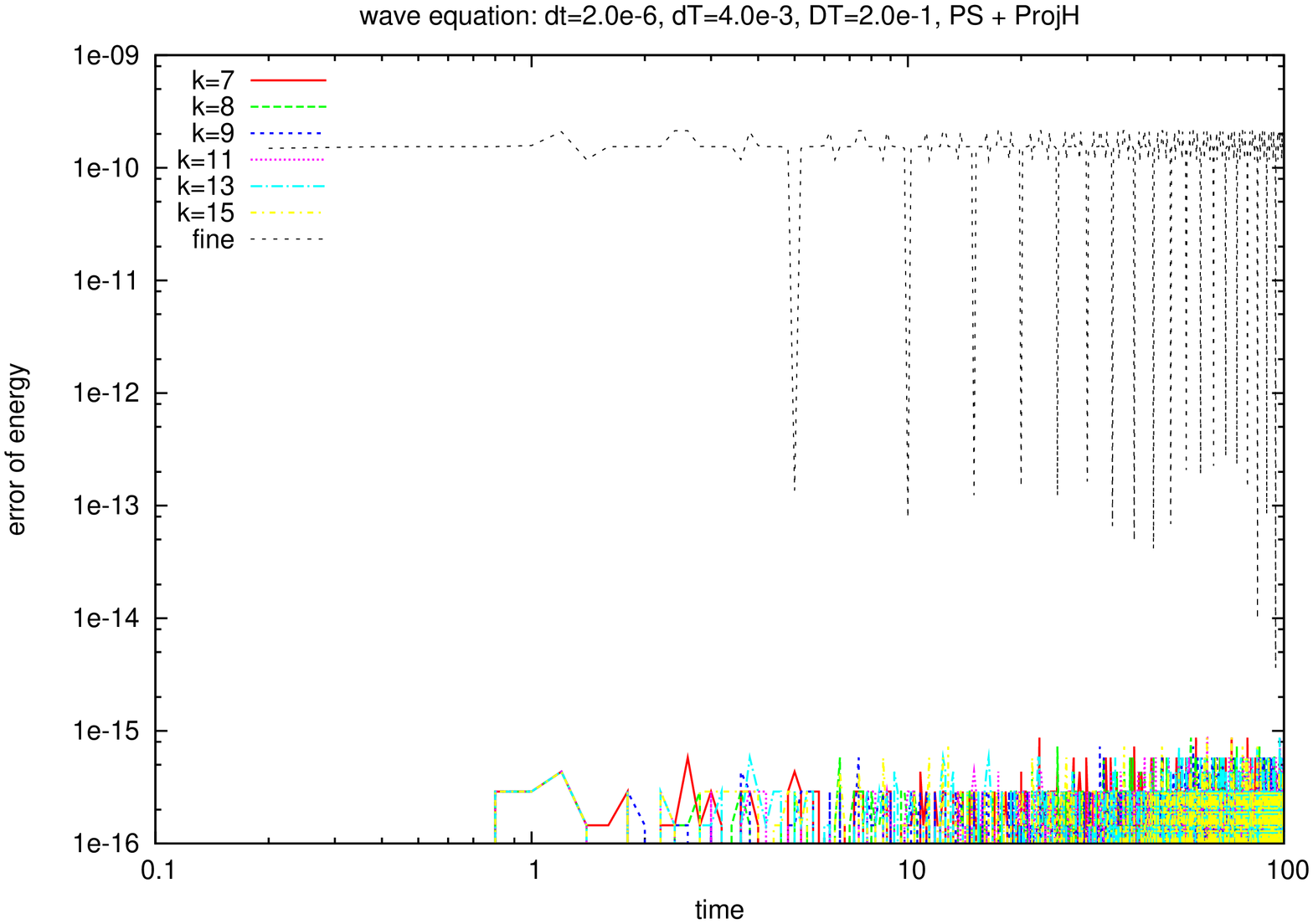}
\caption{Error on energy obtained by parareal method with projection
to the manifold ${\cal M}$ ($\delta t = 2. 10^{-6}, dT=4. 10^{-3},
\Delta T=2. 10^{-1} $)}\label{fig-wave-case1-2-PS-ProjH-E-7-15}
\end{center}
\end{figure}

In the results above, when doing the projection, we view all the
Fourier modes as a whole system, and treat them altogether, with a
unique Lagrange multiplier. Remember that each Fourier coefficient
individually has an energy that is preserved as is explained in
Remark \ref{remark1}. Instead of projecting on the manifold ${\cal
M}$, where only one invariant is taken into account, we can have a
manifold defined with more invariant properties.
We have
 defined only two invariants based on two groups in the Fourier
coefficients, those corresponding to high wave numbers and those
corresponding to low wave numbers. That is, we partition all the
Fourier modes into two groups,   one contains the modes $|k|\ge
K_0$, the other one contains all the modes $|k|< K_0$, we then
perform a projection for each group separately  to make the energy
of each group conserved.

The motivation comes from the fact that,  when
we discretize the solution by Fourier series, the contribution of
the high frequency modes is much smaller than that of the low
frequency modes. If we use the same Lagrange multiplier $\lambda$ to
correct all the modes, the contribution of the high frequency
modes may be to magnify them and pollute the solution.

So, let us introduce two  Lagrange multipliers $\lambda_1$
and $\lambda_2$, one is for the high frequency modes, the other is
for the low frequency modes only.

\begin{figure}[!htbp]
\begin{center}
\includegraphics[width = 8cm, angle=-90]
{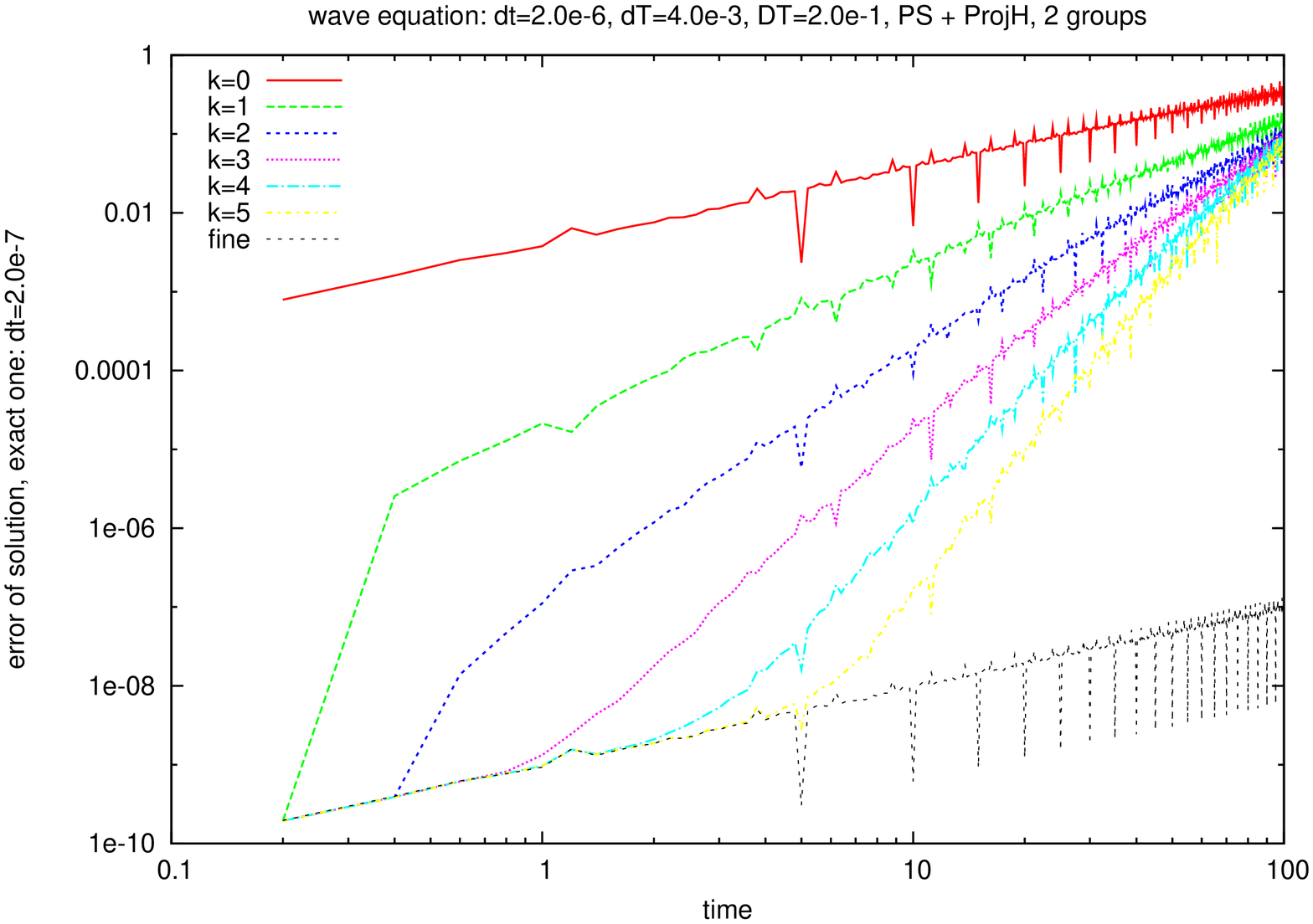}
\caption{Error on the solution obtained by the parareal with
projection on two groups ($\delta t = 2. 10^{-6}, dT=4. 10^{-3},
\Delta T=2. 10^{-1}
$)}\label{fig-wave-case1-0-PS-ProjH-2Group-1-3-P-0-5}
\end{center}
\end{figure}

\begin{figure}[!htbp]
\begin{center}
\includegraphics[width = 8cm, angle=-90]
{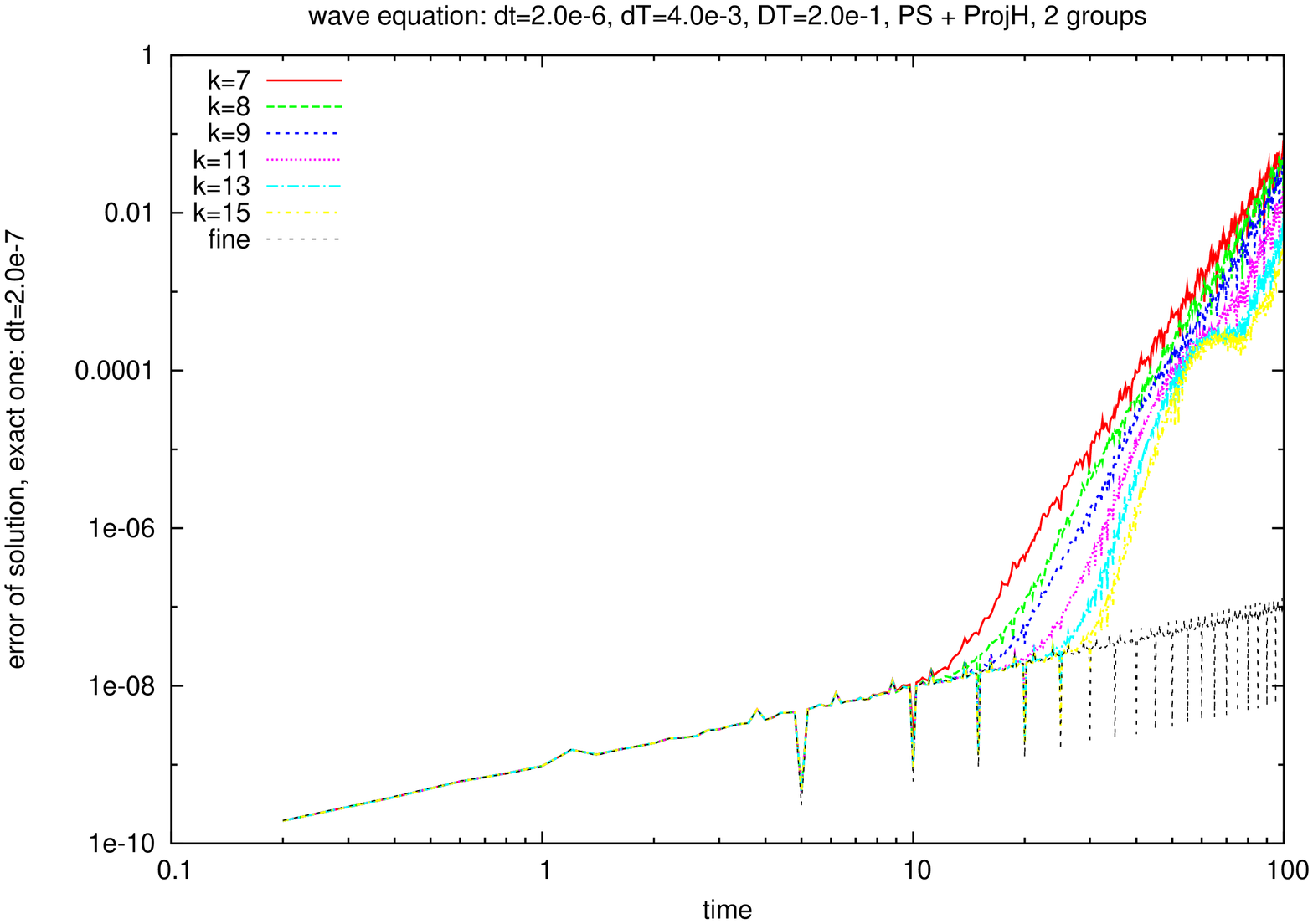}
\caption{Error on the solution obtained by the parareal with
projection on two groups ($\delta t = 2. 10^{-6}, dT=4. 10^{-3},
\Delta T=2. 10^{-1}
$)}\label{fig-wave-case1-0-PS-ProjH-2Group-1-3-P-7-15}
\end{center}
\end{figure}

\begin{figure}[!htbp]
\begin{center}
\includegraphics[width = 8cm, angle=-90]
{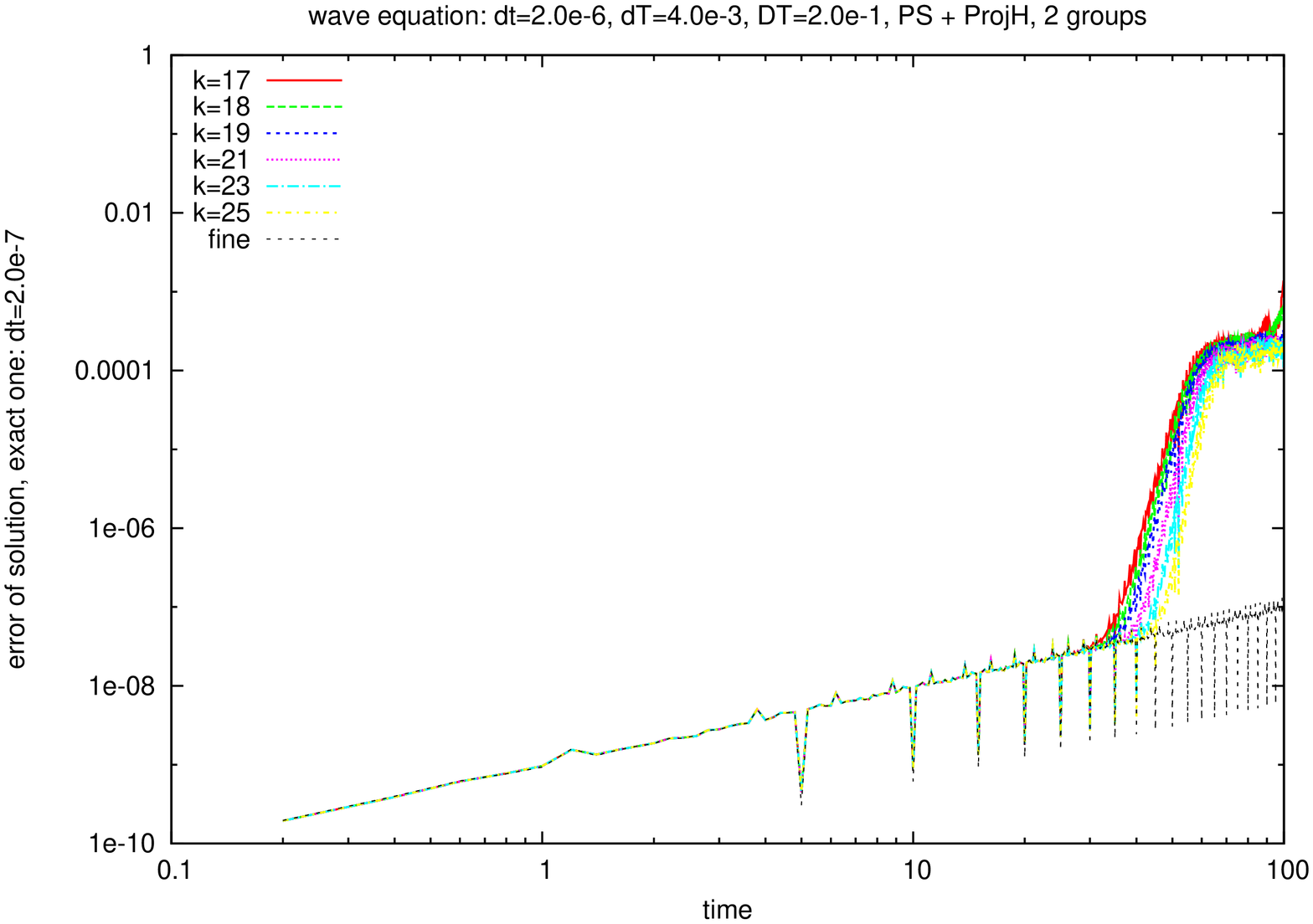}
\caption{Error on the solution obtained by the parareal with
projection on two groups ($\delta t = 2. 10^{-6}, dT=4. 10^{-3},
\Delta T=2. 10^{-1}
$)}\label{fig-wave-case1-0-PS-ProjH-2Group-1-3-P-17-25}
\end{center}
\end{figure}

The method is not more complex to implement. The numerical simulation still with the same $N=30$ is based on $K_0=20$; the results are much improved as can be seen on figures \ref{fig-wave-case1-0-PS-ProjH-2Group-1-3-P-0-5}, \ref{fig-wave-case1-0-PS-ProjH-2Group-1-3-P-7-15} and \ref{fig-wave-case1-0-PS-ProjH-2Group-1-3-P-17-25}.

\begin{figure}[!htbp]
\begin{center}
\includegraphics[width = 8cm, angle=-90]
{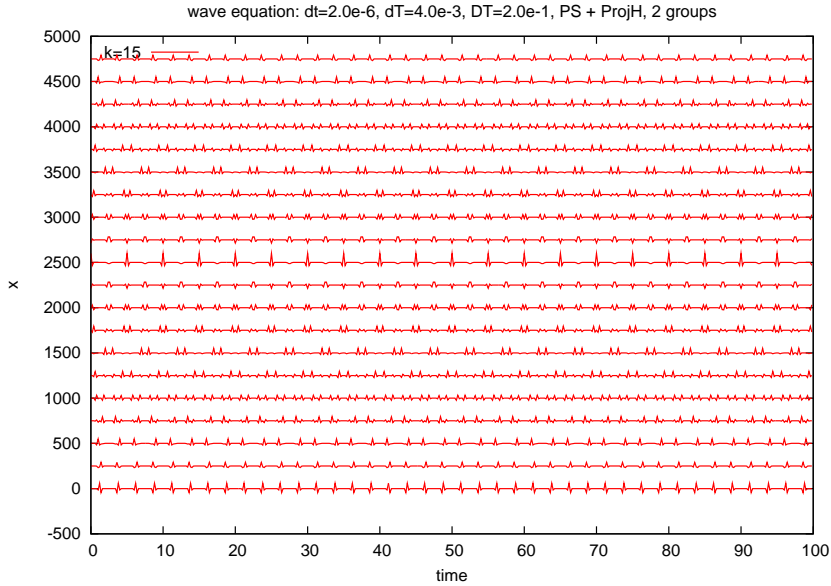}
\caption{Vibration  at different positions obtained by parareal
method with projection on two groups ($\delta t = 2. 10^{-6}, dT=4.
10^{-3}, \Delta T=2. 10^{-1}
$)}\label{fig-wave-case1-0-PS-ProjH-2Group-1-3-vibr-15}
\end{center}
\end{figure}

The solution itself, after 15 iterations, is in excellent accordance with the very fine approximation as can be seen on figure \ref{fig-wave-case1-0-PS-ProjH-2Group-1-3-vibr-15}.

For the linear second order wave equation, the preservation of the invariant quantities based on the decomposition of the coefficients representing the solution into two groups: the low modes and high frequency modes, is thus a sufficiently mild ingredient to improve the behavior of the parareal in time algorithm. In the next section we shall extend this variant of the parareal in time algorithm to a nonlinear PDE.

Note that, if when the method converges, i.e. $\forall n=1,..,N$,  $u_n^k\longrightarrow u_n^\infty$, this limit solution satisfies thanks to (\ref{proj_para})
$$u_{n+1}^{\infty} = \pi_{\cal M}(\mathcal{F}_{\Delta T} (u_{n}^{\infty}))$$
that, here, due to the fact that the propagator $\mathcal{F}_{\Delta T}$ is symplectic and thus preserves the Hamiltonian, reads
$$u_{n+1}^{\infty} = \mathcal{F}_{\Delta T} (u_{n}^{\infty})$$
that proves that the method converges towards the fine solution, as expected.

\section{Burgers equation}

\subsection{The total Fourier/time splitting discretization}

Let us consider in this section the following periodic viscous Burgers' equation in one dimension
\begin{eqnarray} \label{burger-eqaution-1d}
&&\frac{\partial u}{\partial t} - \nu \frac{\partial^2 u}{\partial x^2} + \frac{1}{2} \frac{\partial (u^2)}{\partial x} = 0, \\
&& u(x,0) = u_0(x), \nonumber\\
&& u(x+\mathbb{T},t) = u(x, t), \nonumber
\end{eqnarray}
The viscosity $\nu>0$ is planned to be very small so that the solution will exhibit very sharp gradients.

Taking into account that the periodic frame of the problem allows to use Fourier spectral approximation, we thus consider a semi discretization based on the Galerkin Fourier approximation in $S_N$
\begin{eqnarray} \label{weak-form-burgers-eqaution-1d_N}
(\frac{\partial u_N}{\partial t}, v_N) + \nu(\frac{\partial u_N}{\partial x}, \frac{\partial v_N}{\partial x}) + \frac{1}{2} (\frac{\partial (u_N)^2}{\partial x}, v_N) = 0, \forall v_N
\in  S_N.
\end{eqnarray}
We then choose to use a  second order symmetrized splitting method based on the splitting between the diffusion effects (treated in the Fourier space exactly) and the nonlinear convection effects (treated in the physical space). The only difference between the
coarse propagator $\mathcal{G}$ and fine propagator $\mathcal{F}$ are relative to the  time steps
$dT$ and fine time step $\delta t$, respectively. Denoting by $\delta$ either $dT$ or $\delta t$, the splitting schemes at each iteration proceeds as follows:

\begin{enumerate}
\item Solve the following problem in Fourier space in the interval  $[\tau_n=n\delta, \tau_{n+\frac{1}{2}} = \tau_n+\frac{\delta}{2}]$
\begin{eqnarray}\label{burgers-diff}
(\frac{\partial u_N}{\partial t}, v_N) + \nu (\frac{\partial u_N}{\partial x}, \frac{\partial v_N}{\partial x})  &=& 0,  \forall v_N
\in  S_N\\
u_N(\tau_n) &=& u_{N,n},
\end{eqnarray}
and denote $\tilde{u}_{N,n+\frac{1}{2}} = u_N(\tau_{n+\frac{1}{2}})$. Note that, as in the discretization of the wave equation, by choosing iteratively $v_N=e^{ik\frac{2\pi}{\mathbb{T}}x}$, $k \in [-N, N]$, this results in a system of differential equations in the Fourier coefficients.
\item Solve the following problem  in $[\tau_n, \tau_{n+1}]$
\begin{eqnarray}\label{burgers-convection}
(\frac{\partial u_N}{\partial t},v_N)  + \frac{1}{2} (\frac{\partial (u_N)^2}{\partial x}, v_N) &=& 0,  \forall v_N
\in  S_N\\
u_N(\tau_n) &=& \tilde{u}_{N,n+\frac{1}{2}},
\end{eqnarray}
and denote $u_N(\tau_{n+1})$ as $u_{N,n+\frac{1}{2}}$. Note that here, the solution procedure cannot be exact since the problem is non linear: the discretization in time is based on a $3^{rd}$ order Runge Kutta method and an exact evaluation of the Fourier coefficients of  $\frac{\partial (u_N)^2}{\partial x}$ by a Discrete Fourier Method based on $S_{2N}$.
\item Solve the following problem in Fourier space  in $[\tau_{n+\frac{1}{2}}, \tau_{n+1}]$
\begin{eqnarray}\label{burgers-diff-2}
(\frac{\partial u_N}{\partial t}, v_N) + \nu (\frac{\partial u_N}{\partial x}, \frac{\partial v_N}{\partial x})   &=& 0,\forall v_N
\in  S_N\\
u_N(\tau_{n+\frac{1}{2}}) &=& u_{N,n+\frac{1}{2}}.
\end{eqnarray}
\end{enumerate}

\subsection{Results for the plain parareal in time algorithm}
 
We have performed 3 tests for the plain algorithm. First we note that, provided that the viscosity is large enough, the parareal in time algorithm is working fine. Indeed in case referred later as
{\bf Case 1}, we have chosen

{\bf Case 1:}
\begin{enumerate}
\item $\mu = 1e-2$,
\item $N_c = 2^{8}, N_f = 2^{8}$,
\item{$T=2, \Delta T = 2. 10^{-2}, dT = 1. 10^{-2}, \delta t = 1. 10^{-4}$}
\end{enumerate}

\begin{figure}[!htbp]
\begin{center}
\includegraphics[width = 8cm, angle=-90]
{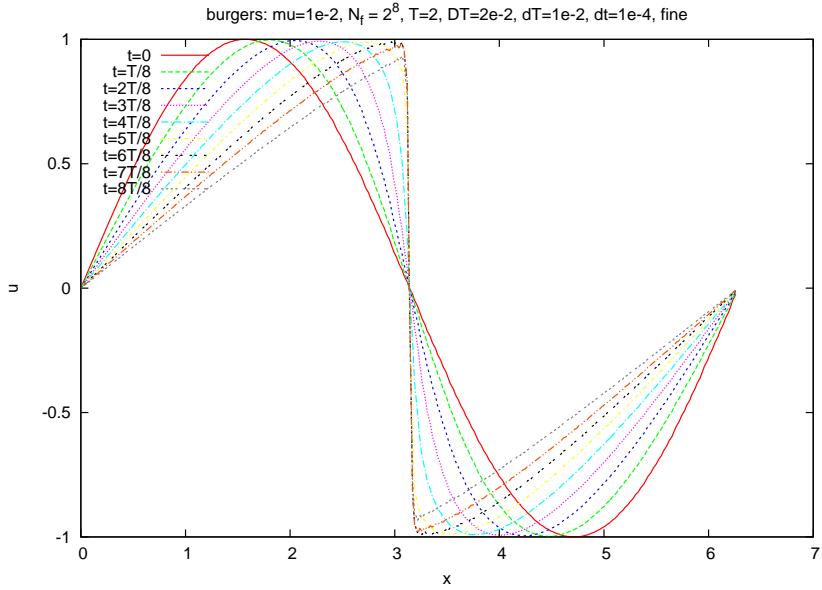}
\caption {Case 1: solution obtained by the parareal method($\delta t
= 1e-4, dT=1e-2, DT=2e-2$)}\label{fig-burgers-case1-solution-fine}
\end{center}
\end{figure}

\begin{figure}[!htbp]
\begin{center}
\includegraphics[width = 8cm, angle=-90]
{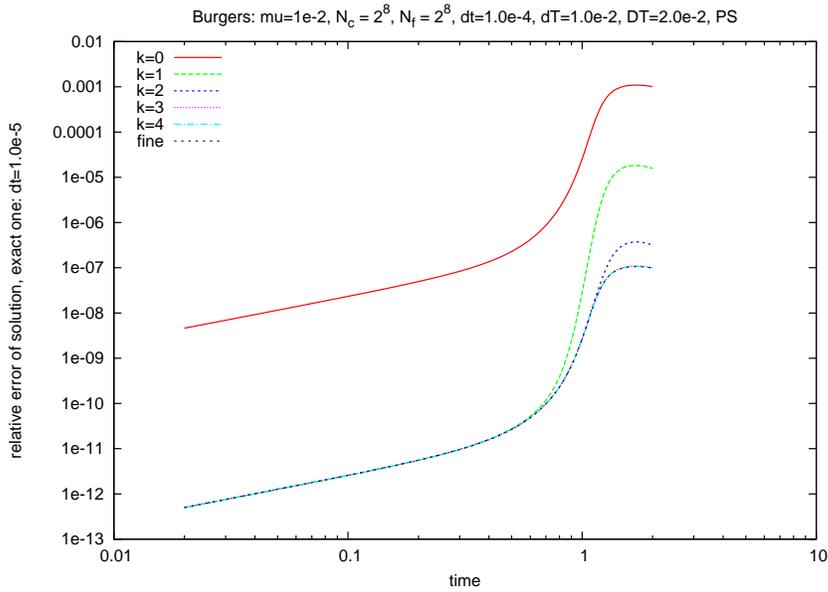}
\caption {Case 1: relative error on the solution obtained by the
parareal method($\delta t = 1e-4, dT=1e-2,
DT=2e-2$)}\label{fig-burgers-case1-PS-SR-0-5}
\end{center}
\end{figure}

The evolution is presented on figure
\ref{fig-burgers-case1-solution-fine}, the relative error in the
solution with respect to a  solution computed with a much smaller time
step is presented on figure \ref{fig-burgers-case1-PS-SR-0-5}
showing that after 3 parareal iterations, the same convergence as
for the fine sequential simulation is achieved.

\bigskip

\bigskip

For the second test, we have diminished the value of the viscosity. This case, referred to as {\bf Case 2}, corresponds to the choices

{\bf Case 2:}
\begin{enumerate}
\item $\mu = 1e-3$;
\item $N_c = 2^{9}, N_f = 2^{9}$,
\item{$T=2, \Delta T = 2. 10^{-2}, dT = 1. 10^{-2}, \delta t = 1. 10^{-5}$}
\end{enumerate}

\begin{figure}[!htbp]
\begin{center}
\includegraphics[width = 8cm, angle=-90]
{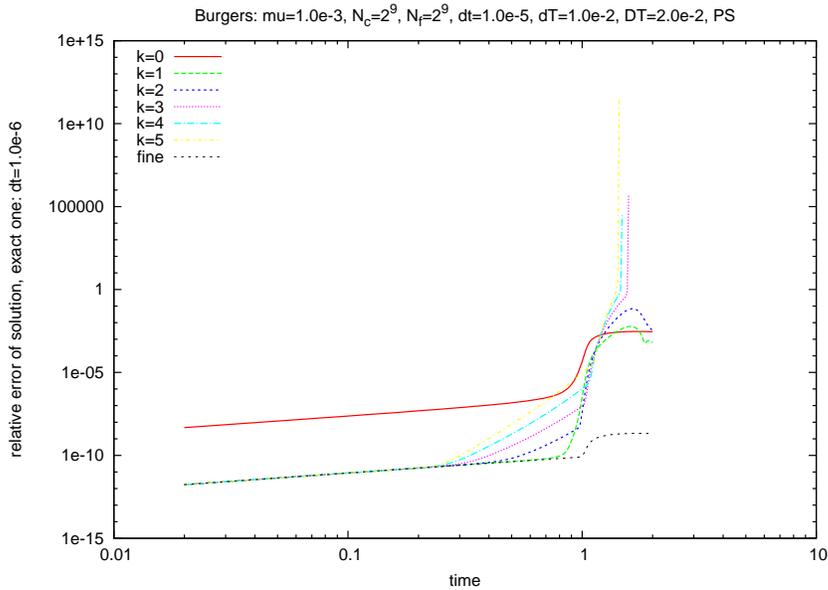}
\caption {Case 2-2: relative error of the solution obtained by the
parareal method($\delta t = 1e-5, dT=1e-2, \Delta
T=2e-2$)}\label{fig-burgers-case2-2-PS-SR-0-5}
\end{center}
\end{figure}

For this case, the plain parareal method is not stable, as can be seen on the figure \ref{fig-burgers-case2-2-PS-SR-0-5}, after 3 iterations instabilities occur soon after time 1 where the solution of the inviscid Burgers equation becomes discontinuous.

It should be noted that, as explained in \cite{MRS}, that using a coarser spacial representation helps in the stabilization of the process. Indeed, {\bf Case 3} below provides a stable parareal approximation

{\bf Case 3:}
\begin{enumerate}
\item $\mu = 1e-3$;
\item $N_c = 2^{8}, N_f = 2^{9}$,
\item{$T=2, \Delta T = 2. 10^{-2}, dT = 1. 10^{-2}, \delta t = 1. 10^{-5}$}.
\end{enumerate}

\begin{figure}[!htbp]
\begin{center}
\includegraphics[width = 8cm, angle=-90]
{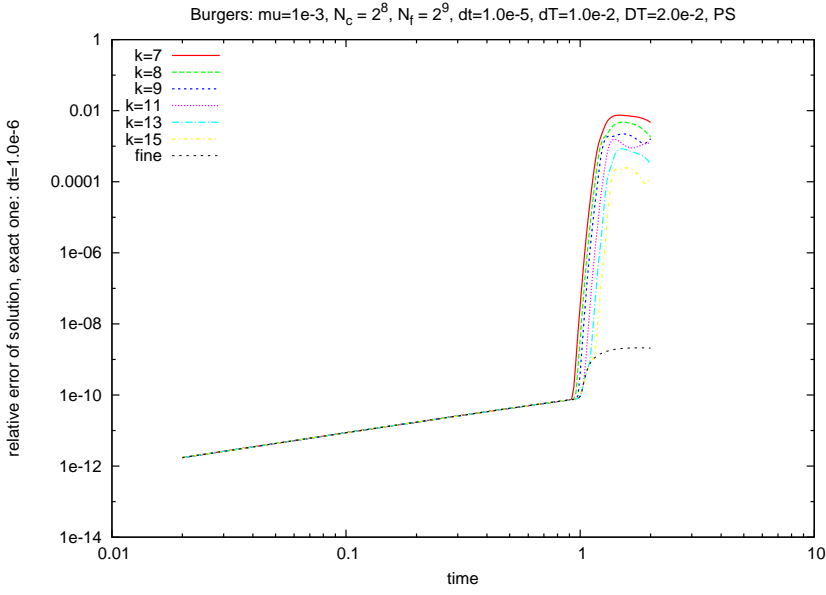}
\caption {Case 2-1: relative error of the solution obtained by the
parareal method($\delta t = 1e-5, dT=1e-2,
DT=2e-2$)}\label{fig-burgers-case2-1-PS-SR-0-5}
\end{center}
\end{figure}

\begin{figure}[!htbp]
\begin{center}
\includegraphics[width = 8cm, angle=-90]
{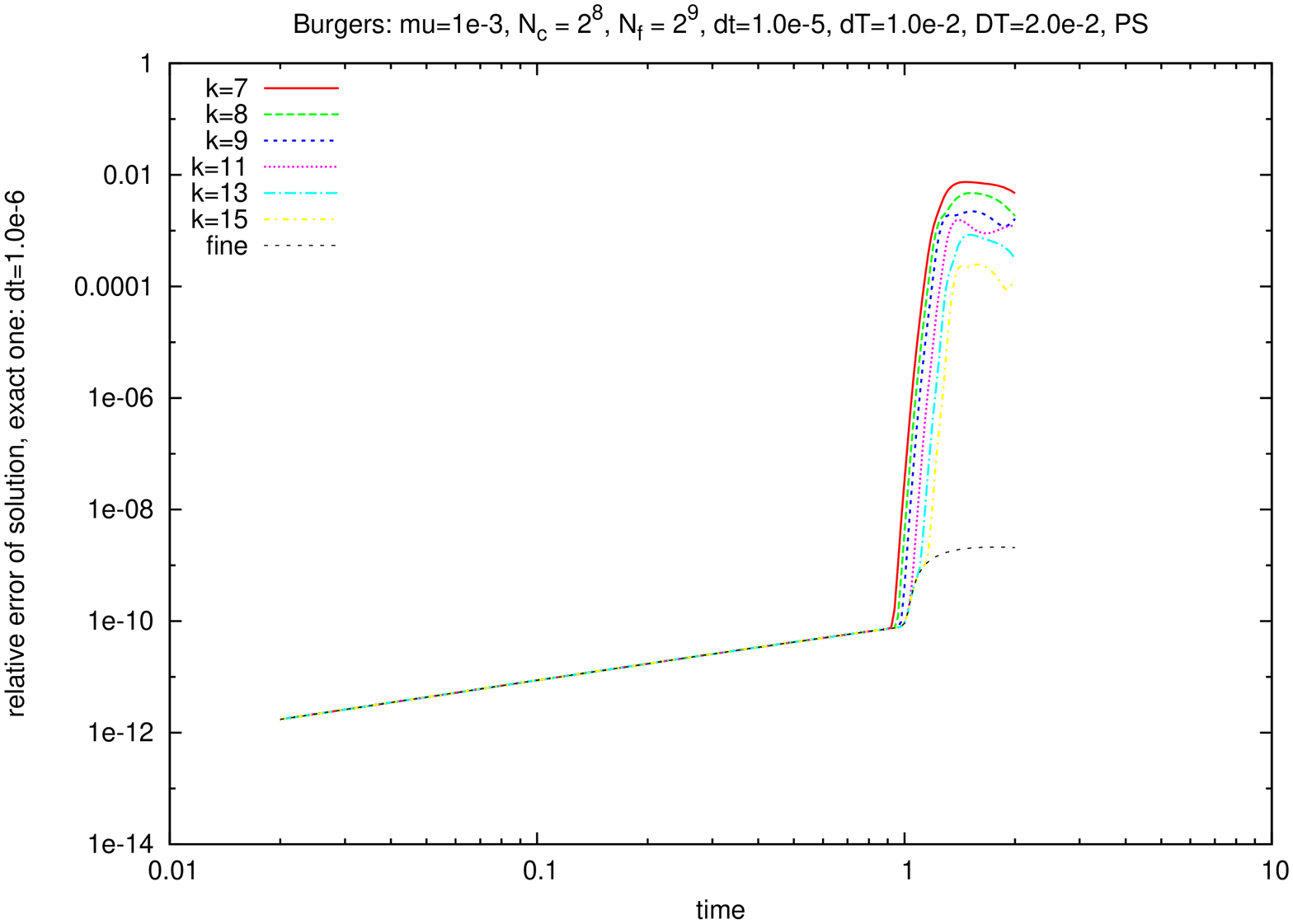}
\caption {Case 2-1: relative error of the solution obtained by the
parareal method($\delta t = 1e-5, dT=1e-2,
DT=2e-2$)}\label{fig-burgers-case2-1-PS-SR-7-15}
\end{center}
\end{figure}

\begin{figure}[!htbp]
\begin{center}
\includegraphics[width = 8cm, angle=-90]
{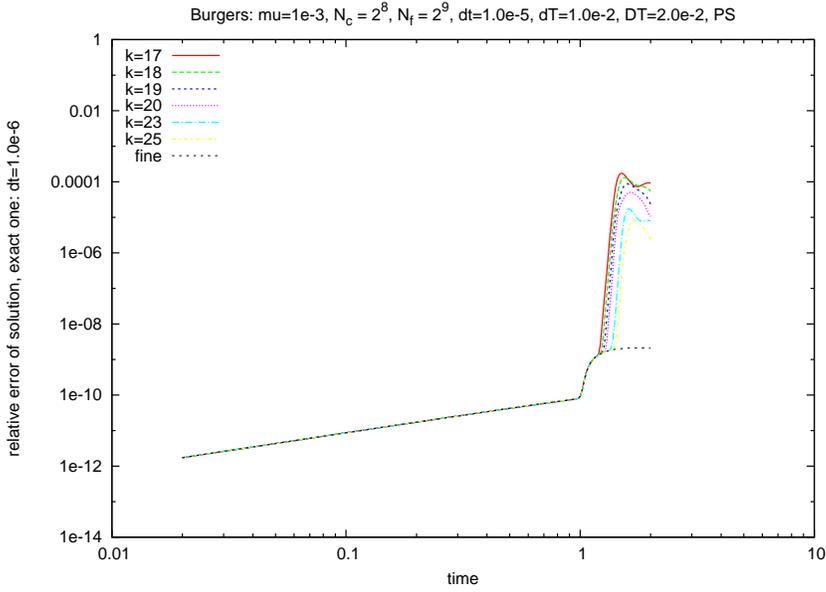}
\caption {Case 2-1: relative error of the solution obtained by the
parareal method($\delta t = 1e-5, dT=1e-2,
DT=2e-2$)}\label{fig-burgers-case2-1-PS-SR-17-25}
\end{center}
\end{figure}

As can be seen on the series of figures \ref{fig-burgers-case2-1-PS-SR-0-5}, \ref{fig-burgers-case2-1-PS-SR-7-15} and \ref{fig-burgers-case2-1-PS-SR-17-25}, the method is stable, converges but very slowly.

In the next sections, we shall improve the two last simulations by extending the method presented in the previous section.

\subsection{Results for the plain parareal in time algorithm with projection alone}

To start with, let us notice that, due to dissipation, there is no quantity that is preserved along a trajectory, nevertheless, by using $u_N$ as a test function in (\ref{weak-form-burgers-eqaution-1d_N}), we derive that
\begin{equation*}
\|u_N(t)\|^2_{L^2(0, \mathbb{T})} + 2 \nu \int_0^t [\frac{\partial u_{N}}{\partial x}]^2 = \|u_N(0)\|^2_{L^2(0, \mathbb{T})},
\end{equation*}
and actually for any two times $T_{n+1}>T_n\ge0$
\begin{equation}
\|u_N(T_{n+1})\|^2_{L^2(0, \mathbb{T})} + 2 \nu \int_{T_n}^{T_{n+1}} [\frac{\partial u_{N}}{\partial x}]^2 = \|u_N(T_n)\|^2_{L^2(0, \mathbb{T})} .
\end{equation}
In the parareal in time algorithm, starting from $u_n^{k+1}$ (remember that this represents the Galerkin approximation that should be denoted as $u_{N,n}^{k+1}$), we expect that
\begin{equation}
\|u_{n+1}^{k+1}\|_{L^2(0, \mathbb{T})} = \|{\cal F}_{\Delta T}[u_{n}^{k+1}]\|_{L^2(0, \mathbb{T})}
\end{equation}
Since we do not want to compute ${\cal F}_{\Delta T}[u_{n}^{k+1}]$ at this stage, we approximate $ \|{\cal F}_{\Delta T}[u_{n}^{k+1}]\|_{L^2(0, \mathbb{T})}$ by
\begin{equation}
 \|{\cal F}_{\Delta T}[u_{n}^{k+1}]\|_{L^2(0, \mathbb{T})} \simeq \frac{ \|{\cal F}_{\Delta T}[u_{n}^{k}]\|_{L^2(0, \mathbb{T})}}{\|u_{n}^{k}\|_{L^2(0, \mathbb{T})}} \|u_{n}^{k+1}\|_{L^2(0, \mathbb{T})}.
\end{equation}

What can thus be done is to consider the manifold ${\cal M}_n^{k+1}$ corresponding to the set of all functions $v$ such that $\|v\|_{L^2(0, \mathbb{T})} = \frac{ \|{\cal F}_{\Delta T}[u_{n}^{k}]\|_{L^2(0, \mathbb{T})}}{\|u_{n}^{k}\|_{L^2(0, \mathbb{T})}} \|u_{n}^{k+1}\|_{L^2(0, \mathbb{T})}
$

Second, we have noticed in the analysis of the wave equation that treating all the frequencies with a unique Lagrange multiplier is not effective enough. We thus partition the frequency space of functions in $S_N$ into two (or more if necessary) groups $G_i$, $i=1,2$ (or $i=1,..,I$) defined by
\begin{equation}
G_i = span\{ e^{i\ell \frac{2\pi}{\mathbb{T}}x}, N_i^-\le |\ell| < N_i^+ \},
\end{equation}
where $L_1^-=0$, $0< L_1^+ = L_2^- -1$, $L_2^+=N+1$
(or $L_1^-=0$, $0< L_1^+ = L_2^- -1$, $L_i^- < L_i^+=L_{i+1}^-+1$, $L_I^+=N+1$ if more groups are around). The $L^2$--projection on group $G_i$ is denoted as $P_i$.
In this case we consider the manifold $[{\cal M}^I]_n^{k+1}$ defined as being the intersection
\begin{equation}
[{\cal M}^I]_n^{k+1} = \cap_{i=1}^I [m^i]_n^{k+1}
\end{equation}
where the group manifolds $ [m^i]_n^{k+1} $ are defined as follows
.\begin{equation}
[m^i]_n^{k+1} = \{ v\in S_N, \| P_i(v)\|_{L^2(0, \mathbb{T})} = \frac{ \|P_i[{\cal F}_{\Delta T}[u_{n}^{k}]]\|_{L^2(0, \mathbb{T})}}{\|u_{n}^{k}\|_{L^2(0, \mathbb{T})}} \|u_{n}^{k+1}\|_{L^2(0, \mathbb{T})}\}
\end{equation}
To any given $\tilde v\in S_N$, we associate the element $\pi_{[{\cal M}^I]_n^{k+1}} (\tilde v)$ defined by the  Standard Projection Method on $[{\cal M}^I]_n^{k+1}$, that involves the determination of $I$ Lagrange multipliers.

With this manifold the parareal scheme with projection reads

\begin{eqnarray} \label{proj_para_burg}\left\{
\begin{array}{rcl}
  u_{n+1}^{0} &=& \mathcal{G}_{\Delta T}(u_n^0) \\
\tilde{u}_{n+1}^{k+1} &=& \mathcal{G}_{\Delta T} (u_{n}^{k+1}) +
\mathcal{F}_{\Delta T} (u_{n}^{k}) - \mathcal{G}_{\Delta T}
(u_{n}^{k}),
\\
u_{n+1}^{k+1} &=& \pi_{[{\cal M}^I]_n^{k+1}}(\tilde{u}_{n+1}^{k+1}).
\end{array}
\right.
\end{eqnarray}

\begin{figure}[!htbp]
\begin{center}
\includegraphics[width = 8cm, angle=-90]
{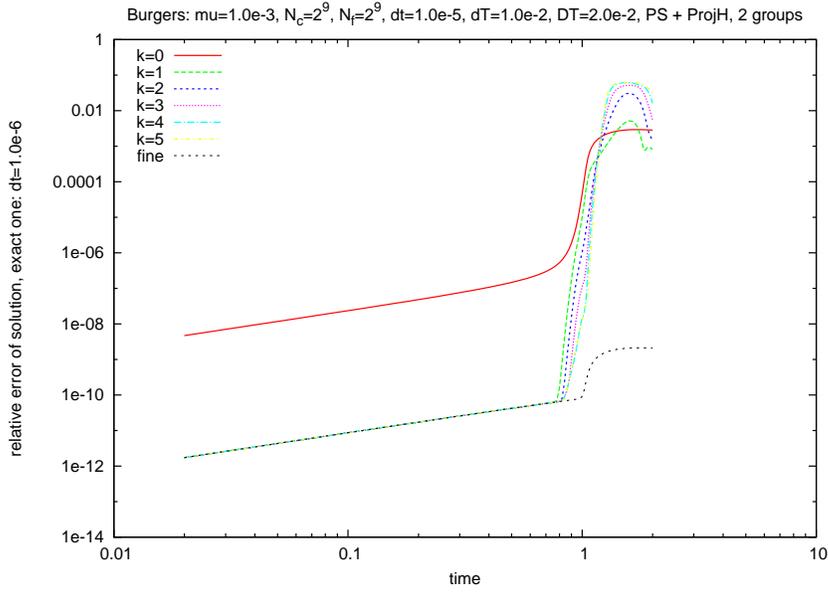}
\caption {Case 2-2: relative error of the solution obtained by the
 method parareal with projection to $2$ groups($\delta t = 1e-5, dT=1e-2,
DT=2e-2$)}\label{fig-burgers-case2-2-PS_ProjH-2group-SR-0-5}
\end{center}
\end{figure}

\begin{figure}[!htbp]
\begin{center}
\includegraphics[width = 8cm, angle=-90]
{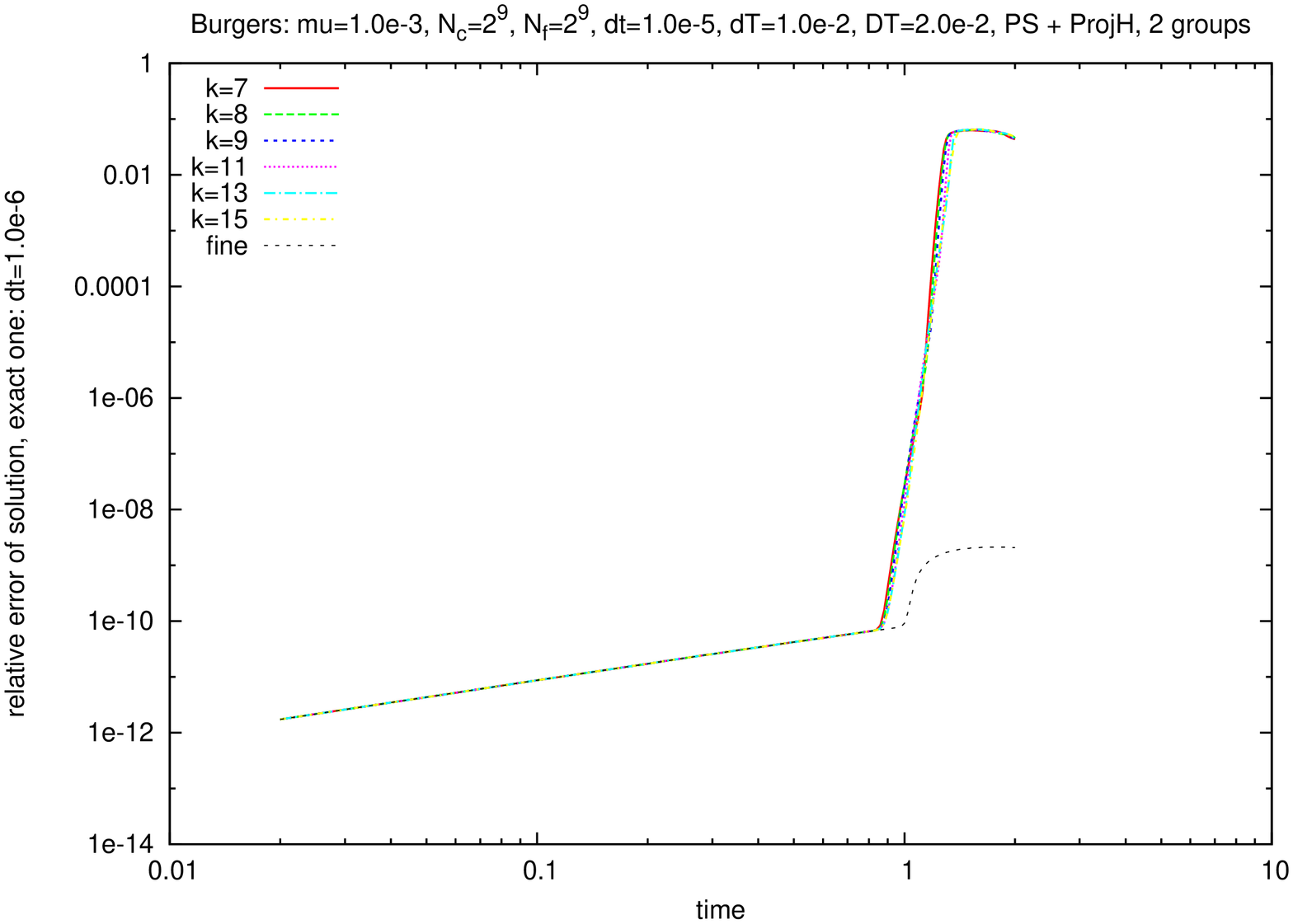}
\caption {Case 2-2: relative error of the solution obtained by the
 method parareal with projection to $2$ groups($\delta t = 1e-5, dT=1e-2,
DT=2e-2$)}\label{fig-burgers-case2-2-PS_ProjH-2group-SR-7-15}
\end{center}
\end{figure}

\begin{figure}[!htbp]
\begin{center}
\includegraphics[width = 8cm, angle=-90]
{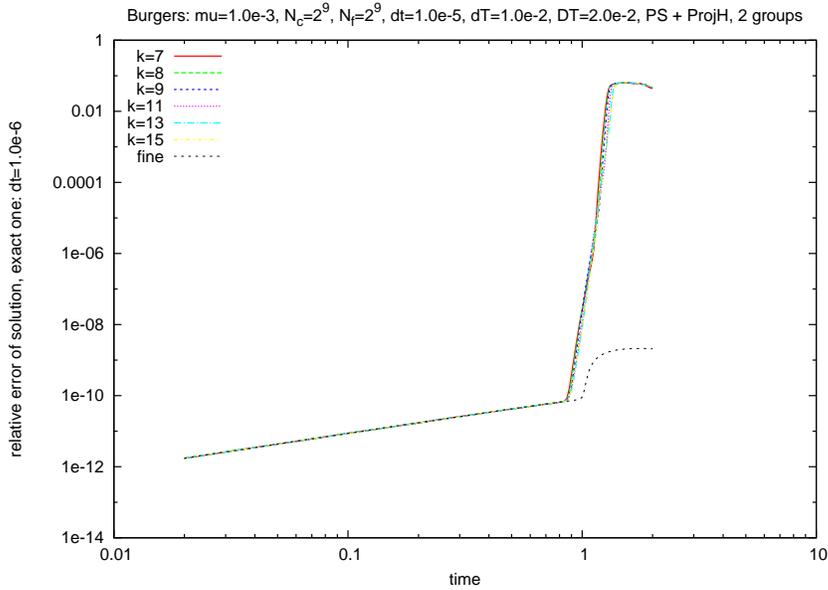}
\caption {Case 2-2: relative error of the solution obtained by the
 method parareal with projection to $2$ groups($\delta t = 1e-5, dT=1e-2,
DT=2e-2$)}\label{fig-burgers-case2-2-PS_ProjH-2group-SR-17-25}
\end{center}
\end{figure}

The results are as follows (see figures \ref{fig-burgers-case2-2-PS_ProjH-2group-SR-0-5}, \ref{fig-burgers-case2-2-PS_ProjH-2group-SR-7-15}, \ref{fig-burgers-case2-2-PS_ProjH-2group-SR-17-25}). We notice that the parareal in time algorithm is stabilized, but it does not converge much. Increasing the number of groups does not improve the behavior of this variant of the parareal in time algorithm.

\begin{remark}
If we would know the exact value of the energy of each groups, the projection method would converge. The problem of course is that due to dissipation, and also to nonlinear effect that mixes up the low and high frequencies, the value of the exact energy in each groups is not known. This is the reason why we get into these difficulties.
\end{remark}

\subsection{Results for the plain parareal in time algorithm with projection and damping}

We notice that, for this non linear problem, the projection allows to stabilize the parareal in time algorithm but a little more is required. We have chosen to damp the high frequencies at each coarse step before proceeding to the next parareal iteration. This means that we introduce a damping function $\Phi$ that, to any element $v_N$ in $S_N$ provides an element $w_N=\Phi(v_N)$ in $S_N$, such that
\begin{eqnarray}
\forall \ell\in G_1,\quad \hat w_\ell &=& \hat v_\ell\\
\forall \ell\in G_i, i=2,..I,\quad \hat w_\ell &=& f(i) \hat v_\ell
\end{eqnarray}
where
$f(i)$ is chosen such that high frequency parts are more and more filtered. In all the following results, we choose $f(i)= \frac{1}{i^2 +1.0}$, except
for $2$ groups, were we have chosen $f(i) = \frac{1}{3 i + 1.0}$.
With the same definition of the multigroup manifold the parareal scheme with projection and damping reads
\begin{eqnarray} \label{proj_para_burg_proj+damp}\left\{
\begin{array}{rcl}
  u_{n+1}^{0} &=& \mathcal{G}_{\Delta T}(u_n^0) \\
\tilde{u}_{n+1}^{k+1} &=& \mathcal{G}_{\Delta T} (\bar u_{n}^{k+1}) +
\mathcal{F}_{\Delta T} (u_{n}^{k}) - \mathcal{G}_{\Delta T}
( \bar u_{n}^{k}), \\
u_{n+1}^{k+1} &=& \pi_{[{\cal M}^I]_n^{k+1}}(\tilde{u}_{n+1}^{k+1})\\
\bar u_{n+1}^{k+1} &=& \Phi( {u}_{n+1}^{k+1}).
\end{array}
\right.
\end{eqnarray}

\begin{figure}[!htbp]
\begin{center}
\includegraphics[width = 8cm, angle=-90]
{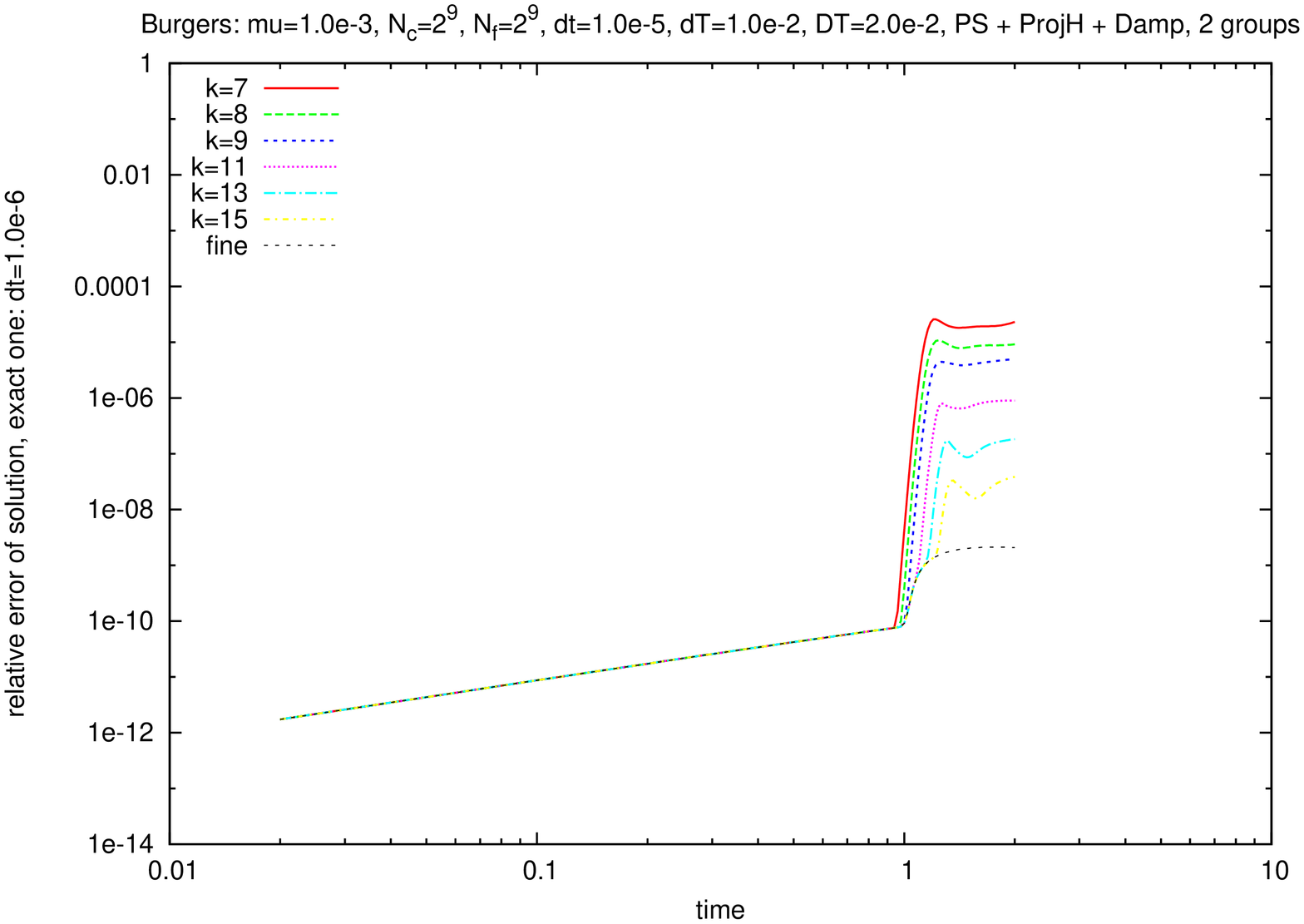}
\caption {Case 2-2: relative error of the solution obtained by the
 method parareal with projection to $2$ groups with damping($\delta t = 1e-5, dT=1e-2,
DT=2e-2$)}\label{fig-burgers-case2-2-PS_ProjH_Damp-2group-SR-7-15}
\end{center}
\end{figure}

The following figures \ref{fig-burgers-case2-2-PS_ProjH_Damp-2group-SR-0-5}, \ref{fig-burgers-case2-2-PS_ProjH_Damp-2group-SR-7-15} and \ref{fig-burgers-case2-2-PS_ProjH_Damp-2group-SR-7-15}
present the results of this improved parareal method. We notice that the parareal method is now converging.

\begin{remark} Note that the strategy used here consisting in adding a little bit of dissipation on the high modes to the parareal scheme goes in the direction suggested in \cite{Bal-2005} and is coherent with the spectral viscosity proposed in \cite{tadmor}. Note also that the solution should be measured after the projection and before the damping. Note finally that the fine propagator acts on the undamped solution so that the dissipativity that is added does not affect the quality of the solution (e.g. as any parareal scheme, after $N^{th}$ iterations, the solution is the same as the one of the plain fine propagator)  but improves the quality of the scheme.
\end{remark}

\begin{figure}[!htbp]
\begin{center}
\includegraphics[width = 8cm, angle=-90]
{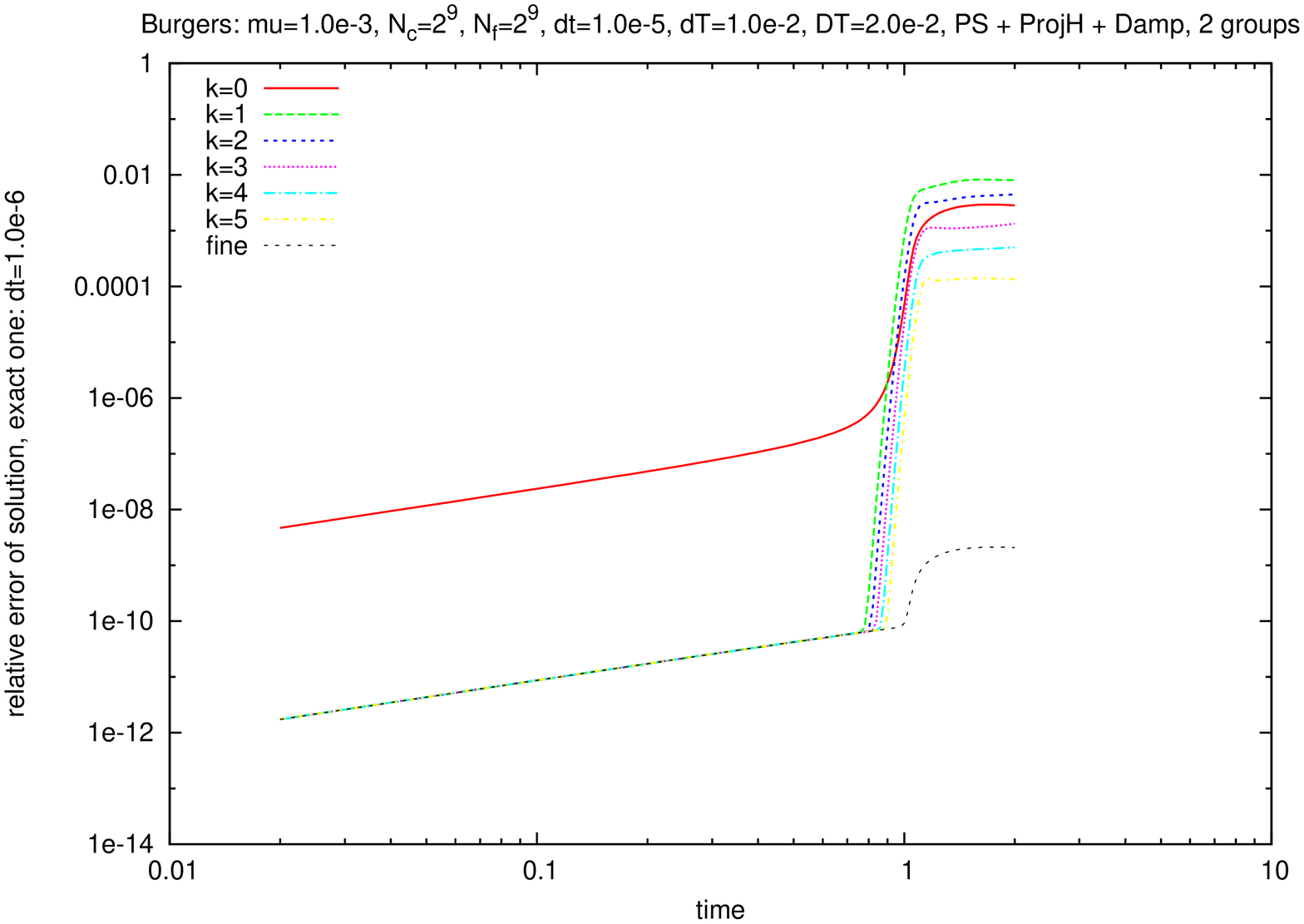}
\caption {Case 2-2: relative error of the solution obtained by the
 method parareal with projection to $2$ groups with damping($\delta t = 1e-5, dT=1e-2,
DT=2e-2$)}\label{fig-burgers-case2-2-PS_ProjH_Damp-2group-SR-0-5}
\end{center}
\end{figure}

\begin{figure}[!htbp]
\begin{center}
\includegraphics[width = 8cm, angle=-90]
{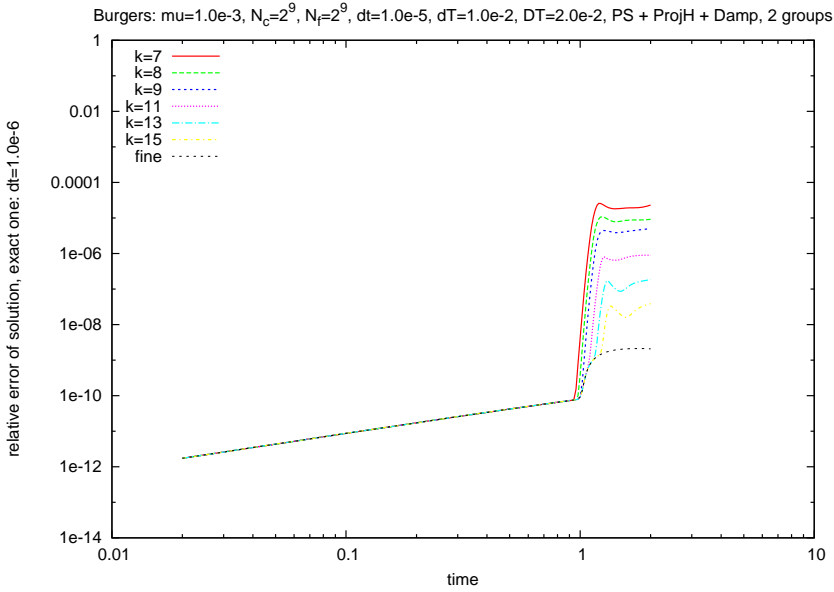}
\caption {Case 2-2: relative error of the solution obtained by the
 method parareal with projection to $2$ groups with damping($\delta t = 1e-5, dT=1e-2,
DT=2e-2$)}\label{fig-burgers-case2-2-PS_ProjH_Damp-2group-SR-17-25}
\end{center}
\end{figure}

\begin{figure}[!htbp]
\begin{center}
\includegraphics[width = 8cm, angle=-90]
{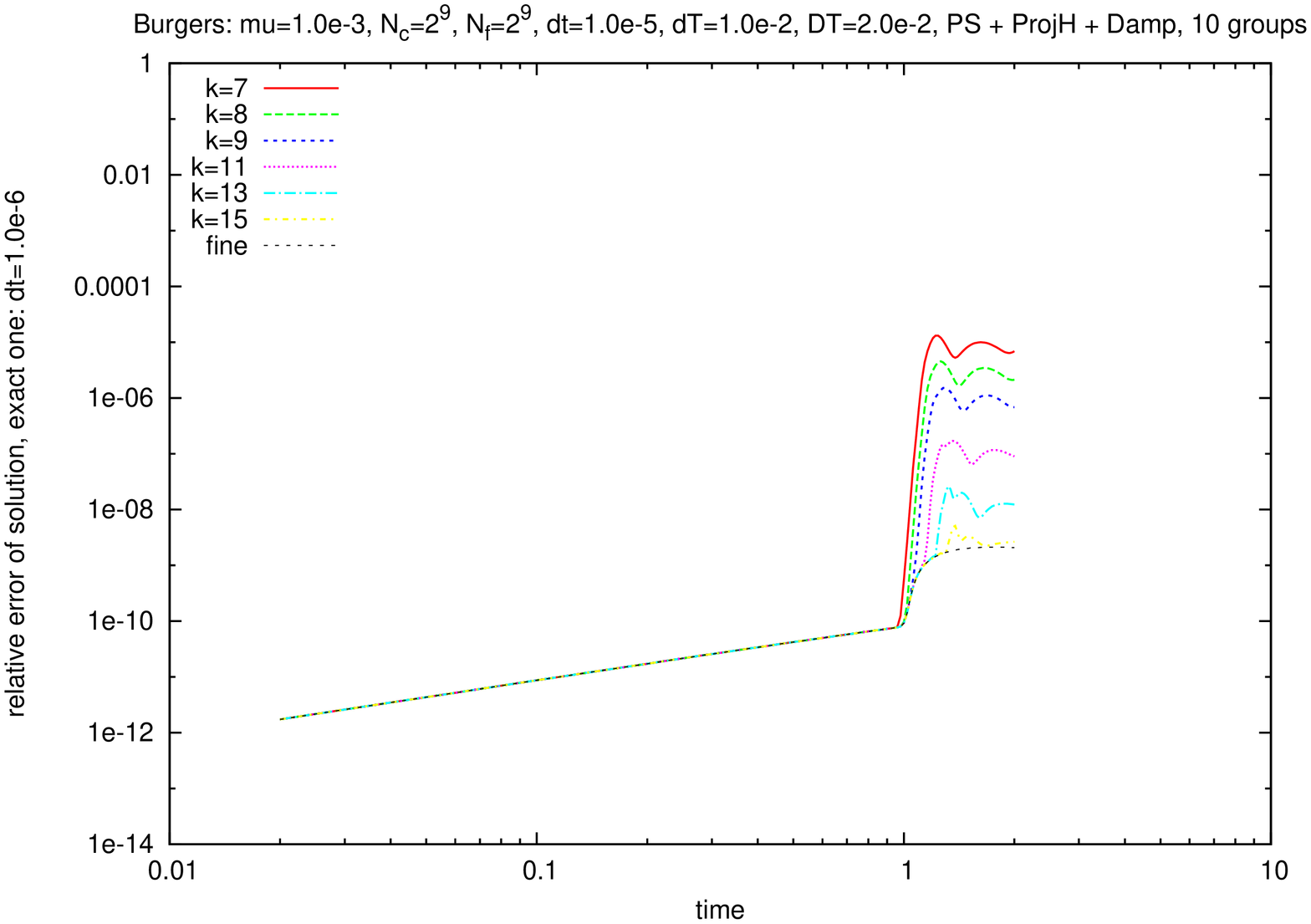}
\caption {Case 2-2: relative error of the solution obtained by the
 method parareal with projection to $10$ groups with damping($\delta t = 1e-5, dT=1e-2,
DT=2e-2$)}\label{fig-burgers-case2-2-PS_ProjH_Damp-10group-SR-7-15}
\end{center}
\end{figure}

The projection with more than two groups performs slightly better, see e.g. figure \ref{fig-burgers-case2-2-PS_ProjH_Damp-10group-SR-7-15} where only the results for iterations between 7 and 15 are presented. It may be interesting to keep this in mind even though the separation in more than two groups is not completely straightforward in the non periodic case where finite difference, finite differences or finite volume techniques will be used.

\begin{figure}[!htbp]
\begin{center}
\includegraphics[width = 8cm, angle=-90]
{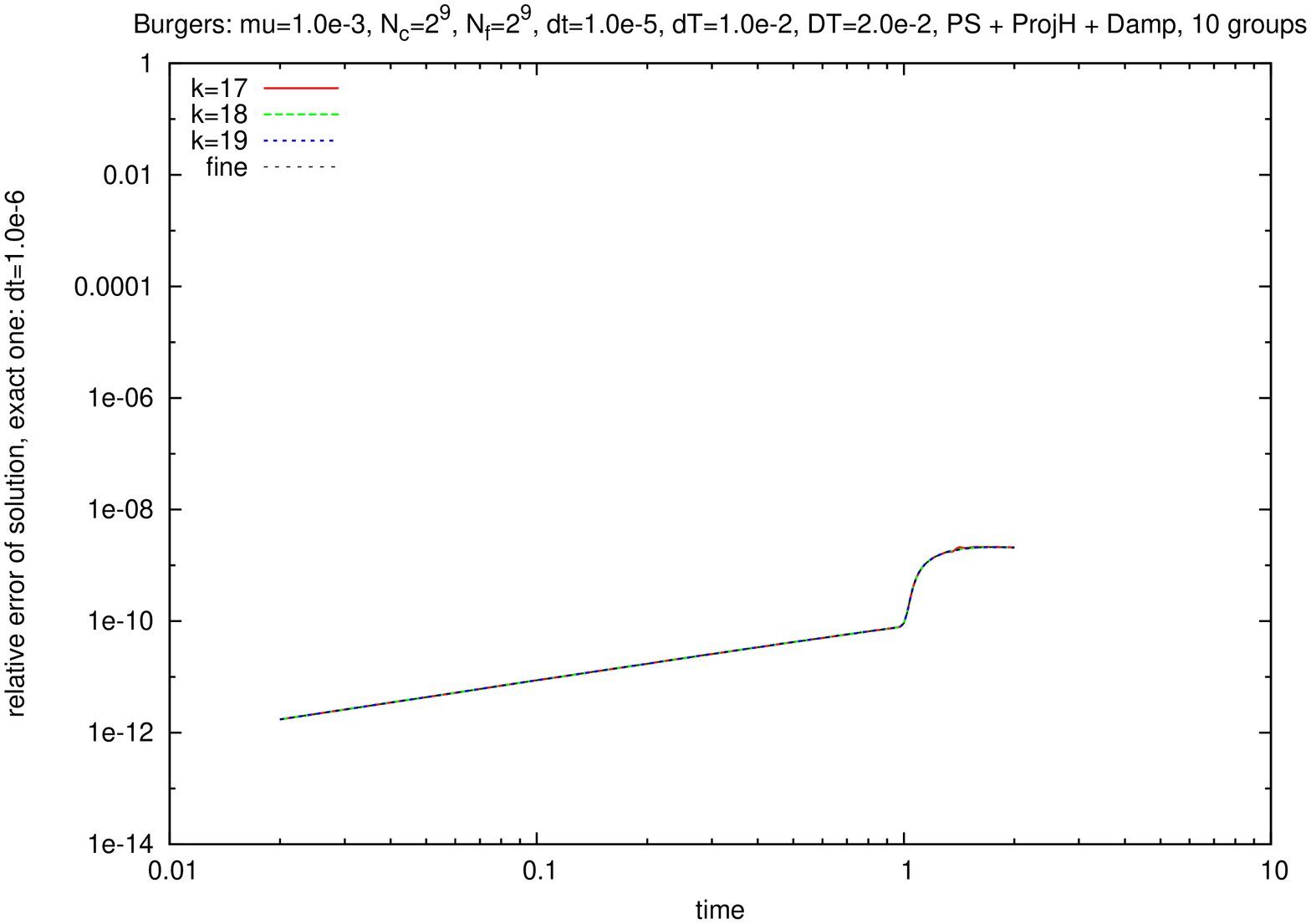}
\caption {Case 2-2: relative error of the solution obtained by the
 method parareal with projection to $10$ groups with damping($\delta t = 1e-5, dT=1e-2,
DT=2e-2$)}\label{fig-burgers-case2-2-PS_ProjH_Damp-10group-SR-17-25}
\end{center}
\end{figure}

\section{Preliminary analysis}

In this section, we provide some simple elements of analysis that prove the stabilty of the previous approximation with the parareal method with projection (and damping).

Let us first assume that the fine propagator is $L^2$ stable in the sense that there exists a constant $\alpha>0$ independent of the discretization parameters such that,
\begin{equation}
\forall v\in L^2(0, \mathbb{T}), \quad \| {\cal F}_{\Delta T} (v)\|_{L^2(0, \mathbb{T})} \le ( 1 +\alpha \Delta T) \| v \|_{L^2(0, \mathbb{T})}
\end{equation}

Note that both for the wave equation and the Burgers' equation, the previous inequality holds with $\alpha = 0$.

Let us start with the case where no dumping is implemented, then, from (\ref{proj_para_burg}) we derive
\begin{eqnarray}
\| u_{n+1}^{k+1} = \pi_{[{\cal M}^I]_n^{k+1}}(\widetilde u_{n+1}^{k+1} )\|^2_{L^2(0, \mathbb{T})} &=& \sum_{i=1}^I \| P_i [u_{n+1}^{k+1} ] \|^2_{L^2(0, \mathbb{T})}\cr
&=& \sum_{i=1}^I \frac{ \|P_i[{\cal F}_{\Delta T}[u_{n}^{k}]]\|^2_{L^2(0, \mathbb{T})}}{\|u_{n}^{k}\|^2_{L^2(0, \mathbb{T})}} \|u_{n}^{k+1}\|^2_{L^2(0, \mathbb{T})}  \cr
&=& \frac{ \|{\cal F}_{\Delta T}[u_{n}^{k}]\|^2_{L^2(0, \mathbb{T})}}{\|u_{n}^{k}\|^2_{L^2(0, \mathbb{T})}} \|u_{n}^{k+1}\|^2_{L^2(0, \mathbb{T})}  \cr
&\le&  ( 1 +\alpha \Delta T)^2 \|u_{n}^{k+1}\|^2_{L^2(0, \mathbb{T})}
\end{eqnarray}
This proves the stability by induction
\begin{eqnarray}
\| u_{n+1}^{k+1} \|_{L^2(0, \mathbb{T})} &\le& ( 1 +\alpha \Delta T)^{n+1} \| u_{0}^{k+1} = u_0 \|_{L^2(0, \mathbb{T})}\cr
&\le& e^{\alpha T} \| u_0 \|_{L^2(0, \mathbb{T})}
\end{eqnarray}

\section{Concluding remarks}

In this paper we have analyzed and proposed a cure to the lack of robustness of the parareal in time algorithm applied to hyperbolic systems or convection diffusion problems with small diffusion. The lack of robustness, leading
sometimes to instabilities of the algorithm before convergence is a consequence of the lack of regularity that the solution develops when time runs. The cure is inherited from the extension of the  parareal in time algorithm for Hamiltonian system and consists in projecting the solution over an energy manifold defined on the fly, and improved iteratively, where the exact solution lies. This simple procedure allows to stabilizes the algorithm that now converges after some iterations. It should be noted that, similarly as for other types of equations, the convergence of the parareal in time algorithm allows to speed up the solution procedure but is not optimal in terms of parallel efficiency since, roughly speaking, a convergence that is achieved after 10 or 15 iterations leads to a parallel efficientcy of $\frac{1}{10}$ or $\frac{1}{15}$ at best, which is far from optimal. It is not the purpose of this paper to provide a full optimal implementation of the approach. Let us just notice that the combination of the parareal in time algorithm with standard domain decomposition methods is a way to improve the efficiency and use more of the processors that are available on current architectures than the one that a plain domain decomposition technique is able to optimally use. We refer to \cite{Guetat-2011} for more on this matter.

\section*{Acknowledgments}
This work has been supported in part by the Agence
Nationale de la Recherche, under the grant ANR-06-CIS6-007 (PITAC).

\end{document}